\documentclass{amsart}
\usepackage[english]{babel}
\usepackage{amsmath}
\usepackage{amsthm}
\usepackage{amscd}
\usepackage{amsfonts}
\usepackage{amssymb}
\usepackage{mathrsfs}
\usepackage{a4wide}
\newcommand{\Ps}{\mathbf{P}^}
\newcommand{\Pry}{\mathcal{P}}

\newcommand{\Z}{\mathbf{Z}}
\newcommand{\M}{\mathcal{M}}
\newcommand{\Sp}{\bar{\mathcal{S}}}
\newcommand{\G}{\mathcal{G}}
\newcommand{\F}{\mathbf{F}}
\newcommand{\C}{\mathbf{C}}

\newcommand{\ra}{\rightarrow}
\newcommand{\str}{\mathcal{O}}
\newcommand{\bircor}{{\dashleftarrow}{\dashrightarrow}}
\renewcommand{\bar}[1]{\overline{#1}}
\renewcommand{\phi}{\varphi}
    \newtheorem{Lem}{Lemma}[section]
    \newtheorem{Prop}[Lem]{Proposition}
    \newtheorem{Thm}[Lem]{Theorem}
    
    \newtheorem{Cor}[Lem]{Corollary}

\theoremstyle{definition}
   \newtheorem{Def}[Lem]{Definition}
    \newtheorem{Exa}[Lem]{Example}
    \newtheorem{Rem}[Lem]{Remark}
    
    \newtheorem{Not}[Lem]{Notation}

    \DeclareMathOperator{\spa}{span}
    \DeclareMathOperator{\mult}{multiplicity}
    \DeclareMathOperator{\sing}{sing}
    \DeclareMathOperator{\Pic}{Pic}
    \DeclareMathOperator{\even}{even}
    \DeclareMathOperator{\odd}{odd}
    
    \DeclareMathOperator{\Jac}{Jac}
    \DeclareMathOperator{\red}{red}
    \DeclareMathOperator{\Div}{Div}
    \DeclareMathOperator{\Aut}{Aut}
\DeclareMathOperator{\Ima}{Im}
\begin{document}
\title{A new model for the theta divisor of the cubic threefold}
\author{Michela Artebani}
\address{Dipartimento di Matematica, Universit\`a di Genova, via Dodecaneso 35, 16146 Genova, Italia}
\email{artebani@dimat.unipv.it}
\author{Remke Kloosterman}
\address{Department of Mathematics and Computer Science, University of Groningen, PO Box 800, 9700 AV  Groningen, The Netherlands}
\email{r.n.kloosterman@math.rug.nl}
\author{Marco Pacini}
\address{Dipartimento di Matematica Guido Castelnuovo, Universit\`a Roma La Sapienza, piazzale Aldo Moro 2, 00185 Roma, Italia}
\email{pacini@mat.uniroma1.it}
\thanks{The authors thank Igor Dolgachev and Alessandro Verra for posing this problem and many useful conversations.
The second author wishes to thank Jaap Top for several useful discussions. 
A part of this research was done during Pragmatic 2003. The authors thank the organizers of this summer school and the University of Catania for giving the possibility to attend this school. The second author thanks EAGER for giving financial support to attend Pragmatic 2003.
}
\date{\today}
\begin{abstract}
In this paper we give a birational model for the theta divisor of the intermediate Jacobian of a generic cubic threefold $X$.
We use the standard realization of $X$ as a conic bundle and a $4-$dimensional family of plane quartics which are totally tangent to the discriminant quintic curve of such a conic bundle structure. The additional data of an even theta characteristic on the curves in the family gives us a model for the theta divisor.
\end{abstract}
\maketitle
\section{Introduction}

Most of the notions mentioned in this introduction are defined in Sections~\ref{spin} (curves), \ref{secDelPezzo} (surfaces) and \ref{oldmodel} (threefolds).

In this paper we give a new birational model for the theta divisor of the intermediate Jacobian of a generic  cubic threefold $X$. 
In \cite[Section 4]{C} a birational model for the theta divisor is given in terms of linear systems of skew cubics on hyperplane sections of $X$. We give a model in terms of even theta characteristics on plane quartics.

Consider the triple $(S,p,D)$, with $S$ a smooth cubic surface, $p$ a point on $S$, not on any line of $S$ and $D$ a double six on $S$. We can associate to such a triple a smooth plane quartic $C$ together with a bitangent $b$ and an Aronhold set $A$ containing $b$. Moreover, we can recover $(S,p,D)$ from the triple $(C,b,A)$. This observation is the main ingredient of our construction.

We try to extend this correspondence as far as possible. The notion of double six can be generalized for a (singular) cubic surface 
with at most isolated $ADE$ singularities.
If $S$ has an elliptic singularity, contains a double line or is reducible, then it seems harder to define  degenerations of double sixes. In some of these cases the possible limit position of the twelve lines giving a double six seems to depend on the degenerating family, hence an intrinsic definition of double six would be impossible.

On the side of plane quartic curves, we need to generalize the notion of bitangent and Aronhold set. The former can be easily defined, while the latter is harder to generalize. Contrary to the smooth case, we need to form generalized Aronhold sets using both generalized bitangents and some components of a curve $\tilde{C}$ associated to $C$. This is enough to give a correspondence between generalized Aronhold sets and generalized double sixes. We show that if $C$ is a stable 
curve, then our generalization coincides with that given by Cornalba in \cite[Section 2]{Corn}. 

Let $X\subset \Ps4$ be a generic cubic threefold and $\ell\subset X$ be a generic line. Let $\tilde{\pi}_\ell : \tilde{X} \ra \Ps2$ be the resolution of the projection from $\ell$. The discriminant curve of the conic bundle $\tilde{\pi}_\ell$ is a smooth plane quintic $Q$. If $E$ is the exceptional divisor of $\tilde{X} \ra X,$ then $\tilde\pi_\ell|_E$ is a finite covering of $\Ps2$ branched over a smooth conic $T.$ Moreover $T\cdot Q=2t$ so that $\theta:=\mathcal O_Q(t)$ is an odd theta characteristic. The line bundle $\theta(-1)$ defines a non-split \'etale double covering $\tilde Q\ra Q$ such that $\tilde Q$ parameterizes the irreducible component of the fibers of $\tilde\pi_\ell$ over $Q.$

Our model for the theta divisor is given by the following:

\begin{Thm}
Let $W$ be the set of quartic curves $C$ such that $d_C=\frac{1}{2}(Q\cdot C)\sim \theta(1)$ and $C$ has at most $ADE$ singularities. Denote with $W^{0}$ the locus of smooth quartics in $W$. 

If $C\in W^{0}$ consider the odd theta characteristic ${\theta}_C=3K_C-d_C$ with associated bitangent $b_C$. Let $B^0=\{(C,b_C)\mid C\in W^0\}$ and $B$ its closure in $W \times {\Ps2}^{*}.$ 

Then 
\[ V_{X,\ell}=\{ (C,b_C,A) \mid (C,b_C) \in B,\ A \mbox{ generalized even theta characteristic on } (C,b_C)   \}\]
is a birational model for $\Theta/\langle -1 \rangle$, the theta divisor of the intermediate Jacobian of $X$ modulo $-1$.
\end{Thm}

The paper is organized as follows:

In Section~\ref{spin} we give some further generalizations of theta characteristics given by Cornalba.

In Section~\ref{secDelPezzo} we discuss some basic properties of Del Pezzo surfaces. 

In Section~\ref{beau} we discuss the correspondence between pairs $(S,\Pi)$ with $S$ a cubic hypersurface in $\Ps{n}$ and $\Pi\subset C$ such that $\Pi\cong \Ps{n-3}$ and pairs $(Q,\theta)$ with $Q$ a curve of degree $n+1$ and $\theta$ an odd theta characteristic.

In Section~\ref{Comb} we give generalizations of odd and even theta characteristics. We use the properties mentioned in Sections~\ref{spin} to~\ref{beau} to show that our definition coincides with the definition of spin structure (from Section~\ref{spin}) in the case of a stable quartic curve. Moreover, we prove the existence of a correspondence between blow-down linear systems on an associated (nodal) Del Pezzo surface and  generalized Aronhold sets.

In Section~\ref{oldmodel} we compare the model of the theta divisor given by Clemens with our model.

In Section~\ref{newmodel} we describe the new model and give some consequences.

In Section~\ref{stablered} we give a connection between the new model and the stable reduction of curves.

\section{Spin curves}\label{spin}
In this section we discuss degenerations of pairs $(C,\theta)$
with $C$ a smooth curve of genus $g$ and $\theta$ a theta
characteristic, using the concept of spin curve.
Cornalba introduced the notion of spin curve and constructed a moduli space $\Sp_g$ of stable spin curves (see \cite{Corn}) with 
a natural morphism $\pi_g: \Sp_g \ra \overline{\M}_g$ of degree $2^{2g}$. 

First, we recall the definition of theta characteristic.
\begin{Def} Let $C$ be a smooth curve of genus $g$. Let $\theta$ be a line bundle on $C$.
 Then $\theta$ is called a {\em theta characteristic} if $\theta\otimes \theta
  \sim K_C$. A theta characteristic $\theta$ is called {\em odd} (resp. {\em even })
  if $h^0(C,\theta)$ is odd (resp. even).
\end{Def}

It is a classical result that a smooth curve of genus $g$ has $2^{2g}$ theta characteristics, of which
$N^{-}_{g}:=2^{g-1}(2^g-1)$ are odd and $N^{+}_{g}:=2^{g-1}(2^g+1)$
 are even.

\begin{Def}
A {\em semi-stable curve} is a reduced, connected curve with only ordinary
double points as singular points such that every smooth rational component contains at least $2$ nodes.

A {\em stable curve} is a semi-stable curve such that every smooth rational component contains at least $3$ nodes.

Let $C$ be a semi-stable curve and $E$ an irreducible component of $C$. Then $E$ is called an {\em exceptional component} of $C$ if $E$ is smooth and rational, such that $\# (\overline{C-E}\cap E)=2.$

A {\em quasi-stable curve} is a semi-stable curve such that any pair of distinct exceptional components is disjoint.\end{Def}

\begin{Rem} If $C$ is a quasi-stable curve then the stable model of $C$ is obtained by contracting every exceptional component. In particular, its
stable model is unique.\end{Rem}

We recall the definition of stable spin curve and explain how to calculate the scheme-structure of the fiber of $\pi_g$ over the points in $\overline{\M}_g-\M_g$. For the latter part we follow \cite[Section 1.3]{CapCas}, to which we refer for the proofs.

\begin{Def}\label{spindef}
A {\em stable spin curve} is a pair $(Y,\theta)$ with $Y$ a
quasi-stable curve (called {\em the support of the spin curve})
and $\theta$ a line bundle on $Y$ such that
\begin{enumerate} \item the restriction of $\theta $ to each exceptional component $E$
is $\mathcal O_{E}(1)$;
\item if we denote $Z:=\bar{Y-\cup E}$, where we take the union over all exceptional components, then
\[(\theta|_{Z})^{\otimes 2}\simeq\omega_{Z},\]
where $\omega_Z$ is the dualizing sheaf.
\end{enumerate}

A stable spin curve is called {\em even} (resp. {\em odd}) if $h^0(Y,\theta)$ is even (resp. odd).

Let $C$ be a stable curve of arithmetic genus $g$. Let $\Sp_g$ be the moduli space of stable spin curves.  The {\em scheme of spin structures on $C$} is the scheme-theoretical fiber of $\pi_g$ over $[C]\in \overline{\M}_g$ and is denoted by $S_C$. Denote with $S_C^{+}$ (resp. $S_C^{-}$) the scheme of even (resp. odd) spin structures on $C$.
\end{Def}

Cornalba (\cite[Lemma 6.3]{Corn}) showed that $\Sp_g$ has two disjoint irreducible components $\Sp_g^+$ and $\Sp_g^-$ corresponding to even and odd spin curves. 

Fix a stable curve $C$.
We characterize all quasi-stable curves appearing as supports in
$S_{C}$. Let $\nu:C^{\nu}\rightarrow C$ be the normalization map,
$B\subset C$ an irreducible component and $B^{\nu}$ the corresponding component
in $C^{\nu}.$
For every subset $\Delta\subset C_{\sing}$ of nodes, set
$D_{B}:=\nu^{-1}(\Delta\cap B)$. Note that $D_B$ is  a divisor on $B^{\nu}$.
\begin{Def} We say that $\Delta$ is even if $\deg D_{B}$ is even for every irreducible component $B$ of $C$.\end{Def}

\begin{Not}
Let $C$ be a reduced nodal curve. 
We denote with $\Gamma_{C}$
the dual graph of $C$, that is the graph whose vertices are the
irreducible components of $C$ and whose edges are the nodes. An
edge connects two vertices if and only if the two corresponding irreducible components intersect in the corresponding node.
If $\Gamma$ is a graph  we  denote with
$b_{1}({\Gamma})$ the first Betti number of $\Gamma$, that is
\[ b_{1}({\Gamma})=\# \{\mbox{ edges }\}- \# \{\mbox{ vertices }\}+\# \{\mbox{ connected components }\}.\]

Moreover, if $C$ is smooth and irreducible then we  denote with $p_{g}(C)$ its geometric genus. 
\end{Not}

The set of all quasi-stable curves $Y$ having $C$ as stable model
is in bijection with the set of subsets of nodes of $C$.
To such a $Y$ we can associate in a unique way the set
$\Delta_{Y}\subset C_{\sing}$ of nodes corresponding to the nodes
of $Y$ not contained in an exceptional component. Conversely, for
every $\Delta\subset C_{\sing}$ there is a unique quasi-stable curve
$Y$ with stable model $C$ and $\Delta_{Y}=\Delta$. 
\begin{Prop}
A quasi-stable curve $Y$ is the support of a spin curve in $S_{C}$
if and only if $\Delta_{Y}$ is even. The number of even subsets
of nodes of $C$ 
is $2^{b_{1}(\Gamma_{C})}.$
\end{Prop}
This is proven in \cite{Corn}.
Fix a quasi-stable curve $Y$ which is the support of a spin curve in $S_{C}$ with $\Delta_Y=\Delta$. Denote with $\nu: Y^{\nu}\rightarrow Y$ the normalization map. 

First we describe which line bundles $\eta$ on $Y^\nu$ are the pullback of a line bundle $\theta$ on $Y$ such that $(Y,\theta)$ is a spin curve. These $\eta$ are characterized by the following properties:
\begin{enumerate}
\item[(1)'] For every component $E$  of $Y^\nu$ such that $\nu(E)$ is an exceptional component of $Y$  we have $\eta|_{E}=\mathcal O_{E}(1)$. 
\item[(2)'] For every non-exceptional component $B$ of $Y$ we have $(\theta^{\nu}|_{B^{\nu}})^{\otimes 2}= K_{B^{\nu}}\otimes \mathcal O(D_{B})$.
\end{enumerate}
From this we deduce that the number of choices for $\eta$ is $2^{2\sum p_{g}(B)}$ where $B$ runs through all the irreducible components of $Y$. 
Let $Z$ as in (2). 
If we fix $\eta$ then the set of all $\theta$ on $Y$
satisfying  relations (1) and (2)  and $\nu^*\theta=\eta$ can be calculated using the exact sequence
\[1\rightarrow (\C^{*})^{b_{1}(\Gamma_{Z})}\rightarrow \Pic(Z) \rightarrow \Pic(
Z^{\nu})\rightarrow 0.\] On every node of $Z$ there are
two compatible gluings, hence we have $2^{b_{1}(\Gamma_{Z})}$ of such
line bundles. The multiplicity of  $\theta$ as a
point in $S_C$ is $2^{b_{1}(\Gamma_{C})-b_{1}(\Gamma_{Z})}$ (see \cite[Section 5]{Corn}).

\begin{Prop} Let $C$ be an irreducible curve of arithmetic genus $g$. Suppose $C$ has $n$ nodes.
Then there are exactly
\[ \left( \begin{array}{c} n \\ k \end{array} \right)  2^{2g-n-k} \] points in
$S_C$ of multiplicity $2^k$. If $k<n$ then half of them are odd
and half of them are even. If $k=n$ then $N^{-}_{g-n}$ of them are odd and $N_{g-n}^{+}$
 are even. \end{Prop}
\begin{proof} Since $b_1(\Gamma_C)=n$, there are $2^n$ quasi-stable
curves which occur as the support of a spin curve in $S_C$. Denote
by $Y$ a resolution of $k$ nodes. Then there are $2^{2g-2n}$ line bundles on $Y^\nu$ satisfying $(1)'$ and $(2)'$. There are $2^{n-k}$ gluing conditions,
hence there are $2^{2(g-n)+(n-k)}$ points in $S_C$ with support $Y$. The
multiplicity of $(Y,\theta)$ in $S_C$ equals $2^k$. If $k<n$ then
 half of the
theta characteristics are odd and half of them are even (see for example \cite[Corollary 2.7]{Har}).  If $k=n$
then the component $Z$ in the support of the spin curve is the normalization of $C$.
The number of odd (resp. even) spin structures with support $Y$ equals the number of odd (resp. even) theta characteristics on $Z$. Note that $Z$ has $N^{-}_{g-n}$ odd and $N^{+}_{g-n}$ even
theta characteristics. The multiplicity of such a spin curve
$(Y,\theta)$ in $S_C$ is $2^n$. There are exactly $ \left(
\begin{array}{c} n
\\ k
\end{array} \right) $ sets of $k$ nodes, which gives the
proposition.
\end{proof}
\begin{Exa}\label{irrexa} If $g=3$ then we have the following results:
\[\begin{array}{lcccc}
 & \mbox{smooth} & \mbox{one node} &\mbox{two nodes} &\mbox{three nodes} \\
\mult 1 \even & 36 & 16 & 8 & 4 \\
\mult 1 \odd  & 28 & 16 & 8 & 4 \\
\mult 2 \even &  - & 10 & 8 & 6 \\
\mult 2 \odd  &  - &  6 & 8 & 6 \\
\mult 4 \even &  - &  - & 3 & 3 \\
\mult 4 \odd  &  - &  - & 1 & 3 \\
\mult 8 \even &  - &  - & - & 1 \\
\mult 8 \odd  &  - &  - & - & 0 \\
\end{array} \]
\end{Exa}

\section{Del Pezzo surfaces}\label{secDelPezzo}
In the sequel we construct several cubic surfaces and double covers  of $\Ps2$ ramified along a reduced quartic. The former are  Del Pezzo surfaces of degree 3, the latter of degree 2. These surfaces are the blow-up of $\Ps2$ in 6, resp. 7 points.
In this section we list some properties of Del Pezzo surfaces and discuss degenerations of Del Pezzo surfaces of degree 2 and 3.
Let $S$ be a smooth Del Pezzo surface of degree $d \in \{2,3\}$. This means that $S$ is the blow-up of $\Ps2$ in $9-d$ distinct points $P_1,\dots P_{9-d}$, such that no three of them lie on a line, and no six lie on a conic. The $(9-d)$-uple  $(P_1,\dots P_{9-d})$, (or, equivalently, the corresponding $9-d$ exceptional curves on $S$)  is called a marking of $S$. 
\begin{Not} \label{mark} To mark the Picard group of $S$ we define the following divisors.
\begin{itemize}
\item Let $L$ be the pre-image of a line in $\Ps2$ not passing through the $P_i$.
\item Let $E_i$ be the exceptional divisor corresponding to $P_i$.
\item Let $L_{i,j}, i\leq j$  be the strict transform of the line connecting $P_i$ and $P_j$. In $\Pic(S)$ we have $L_{i,j}= L-E_i-E_j$.
\item If $d=3$ then let $C_i$ be the strict transform of the conic passing through all the $P_k$ except $P_i$. In $\Pic(S)$ we have $C_i=2L-\sum_{t=1}^6 E_t +E_i$.
\item If $d=2$ then let $C_{i,j}$ be the strict transform of the conic passing through all the $P_k$ except for $P_i$ and $P_j$. In $\Pic(S)$ we have $C_{i,j}=2L-\sum_{t=1}^7 E_t +E_i+E_j$.
\item If $d=2$ then let $D_i$ be the strict transform of the cubic passing through $P_1,\dots, P_7$  with a double point in $P_i$. In $\Pic(S)$ we have $D_i=3L-\sum_{t=1}^7 E_t -2 E_i$.
\end{itemize}
\end{Not}
\begin{Rem} \label{geoint} Except for $L$ all the above listed divisors have self-intersection $-1$. 
\end{Rem}

\begin{Rem}\label{Gei}
If $d=2$ there are 56 smooth rational curves $D$ with $D^2=-1$, we call such a curve an {\em exceptional line}. The morphism $\pi: S\ra | -K_S|\cong \Ps2$ is of degree 2 and the ramification locus of $\pi$ is a quartic curve $C(S)$. The 56 lines of $D$  are the irreducible components of the pre-images of the 28 bitangents of $C$. 
The covering involution $\sigma$ associated to $\pi$ is called the \emph{Geiser involution}. 

We have $\sigma(E_i)=D_i, \sigma(C_{k,l})=L_{k,l}$ and $\sigma(L)=8L-3\sum_{t=1}^7 E_t$  (see \cite[Section VII.4]{Dol}, \cite[Section 7]{GrHa}).

\end{Rem}
\begin{Not}\label{Geidef} Suppose $d=2$. Let $C=C(S) \subset \Ps2$ be the quartic curve  from Remark~\ref{Gei}. Let $D$ be a line on $S$. Then we denote by $b_D$ the bitangent corresponding to $D$, i.e., $b_D= \pi(D)$.
\end{Not}


\begin{Rem}\label{cubic} If $d=3$ then there are 27 smooth rational curves $D$ with $D^2=-1$. They are precisely the 27 lines  on the cubic surface $S\subset | -K_S|  \cong \Ps3$.
 Let $p$ be a point on $S$, not on any of the 27 lines. Consider the projection $\pi_p$ of $S$ with  center $p$. Resolving this map gives a morphism $\pi$ from $\tilde{S}$ to $\Ps2$ where $\tilde{S}$ is the  blow up of $S$ in $p$. The morphism $\pi$ coincides with the anti-canonical map, hence is ramified along a quartic curve $C(S,p)=C(\tilde{S})$. The 28 bitangents of $C(\tilde{S})$ are the image of the 27 lines of $S$ under $\pi_p$, plus the images of the exceptional divisor of the blow-up in $p$.\end{Rem}

\begin{Prop}\label{Cubicblowdown} If $d=3$ then there are exactly $72$ ways of obtaining $S$ as the blow-up of six points in $\Ps2$. Below we list the $72$ associated linear systems $| D |$, together with the 6 curves with self-intersection $-1$ which are blown down by $| D |$.
\[ \begin{array}{ll}
L & \{E_1, E_2, E_3, E_4, E_5, E_6 \} \\
2L-E_l-E_m-E_n&\{E_i, E_j, E_k, L_{l,m}, L_{m,n}, L_{l,n}\} \\
3L-\sum _{t=1}^6 E_t +E_i-E_j&\{E_i, C_i, L_{j,k}, L_{j, l}, L_{j,m},  L_{j,n} \}\\
4L-\sum_{t=1}^6 E_t -E_l-E_m-E_n&\{C_i, C_j, C_k, L_{l,m}, L_{m,n}, L_{l,n}\} \\
5L-2\sum_{t=1}^6 E_t &\{C_1, C_2, C_3,C_4,C_5,C_6\} \\\end{array}\]
In all cases $i,j,k,l,m,n$ are such that $\# \{i,j,k,l,m,n\}=6$.\end{Prop}
\begin{proof} See for example \cite[page 485]{GfHa}.\end{proof}
\begin{Def}\label{DoubleSix} If $d=3$ then a {\em double six} on $S$ is a six-uple of pairs of lines on $S$
\[ ((D_1,D'_1),\dots,(D_6,D'_6))\]
such that $D_i.D_j=D'_i.D'_j=0$, $D_i.D'_j=1$ for $i\neq j$ and $D_i.D'_i=0$.\end{Def}
Using Proposition~\ref{Cubicblowdown} one can show that up to permutation there are  36  double sixes.
\begin{Prop}\label{linsyslist} If $d=2$ then there are exactly 576 ways of obtaining $S$ as the blow-up of seven points on $\Ps2$:
\[ \begin{array}{lll}
L & \{E_1, E_2, E_3, E_4, E_5, E_6, E_7 \} \\
2L-E_m-E_n-E_p&\{E_i, E_j, E_k, E_l, L_{m,n}, L_{n,p}, L_{m,p}\} \\
3L-\sum_{t=1}^7 E_t +E_i+E_j-E_k&\{E_i, E_j, C_{i,j}, L_{k,l}, L_{k, m}, L_{k,n},  L_{k,p} \} \\
4L-\sum_{t=1}^7 E_t +E_i -E_m-E_n-E_p& \{E_i, C_{i,j}, C_{i,k}, C_{i,l}, L_{m,n}, L_{n,p}, L_{m,p}\} \\
5L-2 \sum_{t=1}^7 E_t + 2E_i &\{E_i,C_{i,j}, C_{i,k}, C_{i,l},C_{i,m},C_{i,n},C_{i,p}\} \\
\hline
8L-3\sum_{t=1}^7 E_t & \{D_1, D_2, D_3, D_4, D_5, D_6, D_7 \} \\
7L-2\sum_{t=1}^7 E_t-E_m-E_n-E_p&\{D_i, D_j, D_k, D_l, C_{m,n}, C_{n,p}, C_{m,p}\} \\
6L-2\sum_{t=1}^7 E_t -E_i-E_j+E_k&\{D_i, D_j, L_{i,j}, C_{k,l}, C_{k, m}, C_{k,n},  C_{k,p} \} \\
5L -\sum_{t=1}^7 E_t-2E_i-E_j-E_k-E_l& \{D_i, L_{i,j}, L_{i,k}, L_{i,l}, C_{m,n}, C_{n,p}, C_{m,p}\} \\
4L-\sum_{t=1}^7 E_t -2 E_i&\{D_i,L_{i,j}, L_{i,k}, L_{i,l},L_{i,m},L_{i,n},L_{i,p}\} \\
\end{array}\]
In all cases $i,j,k,l,m,n,p$ are such that $\# \{i,j,k,l,m,n,p\}=7$.
\end{Prop}
\begin{proof} The order of the Weyl group $W(E_7)$ equals the number of ordered $7$-uples of points in $\Ps2$ (up to automorphisms of $\Ps2$) giving the same Del Pezzo surface (cf. \cite[Theorem 3, VI]{Dol}). Since $\# W(E_7)=576 \# S_7$  there are 576 markings up to permutation of the 7 points in $\Ps2$. One can easily show that all the above listed systems give different markings and that their number is 576.
\end{proof}
\begin{Rem}\label{ds} The horizontal line divides the linear systems mentioned in Proposition~\ref{linsyslist}  in two groups. The upper half are those coming from the following construction: let $Y$ be the blow-up of $\Ps2$ at 6 points out of the 7 points $P_i$. List all the linear systems that give a blow-down model $Y\ra\Ps2$. Take then the blow-up of the 7th point, to obtain a blow-down model $S\ra\Ps2$. The lower half is obtained by applying the Geiser involution (see Remark ~\ref{Gei}) to the upper half of the list.

The first 288 items in the list of  7-uples of divisors give rise to 288 sets of 7 bitangents. These are precisely the 288 Aronhold sets associated to $C(S)$. Each Aronhold set determines an even theta characteristic by mapping $\{D_1,\dots,D_7\}$ to $\sum b_{D_i}-3K$. Each even theta characteristic can be obtained from 8 different 7-uples of divisors.
Fix a divisor $D_0$ and consider all linear systems $|D|$ in Proposition~\ref{linsyslist} containing $D_0$, then we find 72 linear systems. These 72 linear systems are divided up in 36 pairs, such that each pair corresponds to an even theta characteristic. Consider the cubic surface $Y'$ obtained by blowing down $D_0$. Then the 72 above linear systems correspond to the 72 linear systems on the cubic surface $Y'$ which give a blow-down model $Y' \ra \Ps2$. The 36 pairs of linear systems correspond to the 36 double six.\end{Rem}

Del Pezzo surfaces can be defined equivalently as nonsingular rational surfaces with ample anti-canonical class $-K_S$. If $d>2$ then the anti-canonical linear system maps the surface isomorphically to a smooth surface of degree $d$ in $\Ps{d}$. As mentioned before, to obtain a Del Pezzo surface of degree $d$, $d\geq 1$, we have to choose $9-d$ distinct points $P_i$ in $\Ps2$ such that no three of them are collinear and no six of them are on a conic. Thus a parameter space for these surfaces  is an open dense subset of $(\Ps2)^{9-d}$. Taking the quotient for the action of $PGL(3)$ gives rise to a moduli space $\M_{DP}(d)$ for marked Del Pezzo surfaces of degree $d\leq 8$ (cf. \cite[Theorem VI.3]{Dol}).

\begin{Def} A {\em nodal Del Pezzo surface} is a smooth surface $S$ with almost ample (i.e., big and nef) anti-canonical class. The {\em degree} of $S$ is $K_S^2$.\end{Def}

Nodal Del Pezzo surfaces can be obtained by taking point sets in $\Ps2$ in `almost general' position, namely,  point sets containing `infinitely near' points, three collinear points or six points on a conic (see \cite[Sections VII.3 \& VII.4]{Dol} for precise conditions on point sets). A parameter space for such surfaces is then given by an open set in the smooth variety $\hat{P}_2^{9-d}$ parameterizing infinitely near point sets in $\Ps2$.

The anti-canonical model of a nodal Del Pezzo surface $S$ of degree $d>2$ is a normal surface in $\Ps{d}$ with rational double points (coming from contracting the $-2$-curves on $S$).  If $d=2$ the anti-canonical map of $S$ can be factored as  the contraction of all the $-2$-curves, composed with a degree 2 morphism to $\Ps2$ ramified along a (possibly singular) quartic curve. 

\begin{Not} We denote with $\mathcal{N}\subset \Pic(S)$ the subgroup generated by  smooth rational curves with self-intersection $-2$.\end{Not}

\begin{Def}\label{Defsch} Let $S$ be a nodal Del Pezzo surface of degree $d$. The {\em scheme of lines} $\mathcal{L}_S$ is the Hilbert scheme of smooth rational curves $D$ on $S$ with $\omega_S\otimes \str_D \cong \str_{\Ps1}(1)$. The {\em multiplicity} of $D$ is
\[\frac{\# \{ \sigma \in \Aut(NS(S)) : \sigma(D)\equiv D \bmod \mathcal{N}\} }{\# \{ \sigma \in \Aut(NS(S)) : \sigma(D)=D \}}. \]

The {\em scheme of blow-down models } $\mathcal{BM}_S$ is the zero-dimensional scheme of blow-down models $|L|: S\ra \Ps2$. The {\em multiplicity} of $|L|$ is
\[\frac{\# \{ \sigma \in \Aut(NS(S)) : \sigma(L)\equiv L \bmod \mathcal{N}\} }{\# \{ \sigma \in \Aut(NS(S)) : \sigma(L)=L \}}. \]

If $d=3$ we can define an involution on $\mathcal{BM}_S$ in the following way. Suppose $|L|$ is blow-down morphism. Then there are six rational curves contracted by $L$ say $F_1$ up to $F_6$. We can form a unique set of divisors  $\{D_1,\dots,D_6\}$, such that the $D_i$ are reduced connected effective divisors,  $D_i.D_j=-\delta_{i,j}$ and $D_i$ contains at least one component which is an exceptional line. Then $\sigma(L)=5L-2\sum D_i \bmod \mathcal{N}$. We define the {\em scheme of double-six $\mathcal{D}_S$} as $\mathcal{BM}_S$ modulo the action of $\sigma$.
\end{Def}

\begin{Rem} To the knowledge of the authors, there is no place in the literature where the schemes $\mathcal{L}_S, \mathcal{BM}_S$ and $\mathcal{D}_S$ are defined. It seems to the authors that this is the most natural definition.

If $S$ were a classical Del Pezzo surface then $\Aut(NS(S))$ would act transitively on all exceptional lines in $S$, and if $D$ is an exceptional line then the divisor $\sigma(D)$ is an exceptional line for every $\sigma \in \Aut(NS(S))$. If $S$ is nodal then there might exists $\sigma\in \Aut(NS(S))$ such that $\sigma(D)$ is not an exceptional line, but one can show that the class $\sigma(D) \mod \mathcal{N}$ contains precisely one exceptional line. Similar reasonings can justify the definition of the other two schemes. 

To obtain a good definition of double six on a nodal Del Pezzo surface of degree 3, one needs to describe this double six in terms of the two associated linear systems. One can construct examples where the position of the $-1$ divisors modulo $\mathcal{N}$ does not determine the double six in the above sense (cf. Example~\ref{exaunique}).
\end{Rem}

\begin{Rem} To show that $\mathcal{D}_S$ is well-defined one has to prove that the set $\{D_1,\dots, D_6\}$ exists and is unique. This is an easy exercise using Dynkin diagrams.\end{Rem}
\begin{Rem} If $d=2$ then the scheme $\mathcal{L}_S$ has length 56 and the scheme $\mathcal{BM}_S$ has length 576. There is a natural action of the Geiser involution on both schemes.

If $d=3$ then the length of $\mathcal{D}_S$ is 36.
\end{Rem}

\section{Determinantal hypersurfaces}\label{beau}
For an overview of this subject we refer to a recent paper of A. Beauville \cite{Beau}. That paper deals with the classical question to determine when an integral hypersurface in $\Ps{n}$ can be written as the zero-set of the determinant of a matrix with homogeneous entries. We are interested in a special case, namely the symmetric determinantal representations of smooth plane curves:
\begin{Prop}[{Beauville, \cite[Proposition 4.2]{Beau}\label{det}}] Let $C$ be a smooth plane curve of degree $d$, defined by the equation $F=0$. Let $\theta$ be a theta characteristic on $C$.
\begin{enumerate}
\item\label{eventh} If $h^0(\theta)=0$ then there exists a minimal resolution of $\theta$ (unique up to isomorphism)
\[ 0 \ra \str_{\Ps2} (-2)^d \stackrel{M}{\ra} \str_{\Ps2} (-1)^d  \ra \theta \ra0,\]
where $M$ is a symmetric matrix of linear forms and $\det(M)=F$.
\item \label{oddth}If $h^0(\theta)=1$ then there exists a minimal resolution of $\theta$ (unique up to isomorphism):
\[ 0 \ra \str_{\Ps2} (-2)^{d-3} \oplus \str_{\Ps2}(-3) \stackrel{M}{\ra} \str_{\Ps2} (-1)^{d-3} \oplus \str_{\Ps2}  \ra \theta \ra0,\]
where $M$ is a symmetric matrix of the form
\[ \left( \begin{array}{cccc}
L_{1,1} & \cdots &\L_{1,d-3} & Q_1 \\
\vdots &&\vdots&\vdots\\
L_{1,d-3}&\cdots&L_{d-3,d-3}&Q_{d-3} \\
Q_1&\cdots&Q_{d-3} & H
\end{array} \right) \]
with $L_{i,j}$ linear forms, $Q_k$ quadratic forms, $H$ a cubic and $\det(M)=F$.
\end{enumerate}
Conversely, the cokernel of the map defined by a matrix $M$ of the form in (\ref{eventh}) (resp. (\ref{oddth})), such that $\det(M)\neq 0$, gives rise to a theta characteristic on the curve given by $\det(M)=0$ with $h^0(\theta)=0$ (resp. $h^0(\theta)=1$.) \end{Prop}
If $d=5$ then the determinantal representation of $C$ through $\theta$ can be endowed with a geometric interpretation.

Fix a couple $(Q,\theta)$ where $Q$ is a smooth plane quintic and $\theta$ a theta characteristic on $Q$ with $h^0(\theta)=1$. The associated matrix $M$ has the form
\[\left( \begin{array}{ccc} L_{1,1}& L_{1,2} & Q_1\\
 L_{1,2} &L_{2,2}&Q_2\\
Q_1&Q_2 &H \end{array} \right).\]
Define the cubic threefold $X=X(C,\theta)$ in $\Ps4$ as the zero-set of
\begin{equation*}
\sum_{i,j\in\{1,2\}}{u_iu_{j}L_{ij}(x_0,x_1,x_2)}+\sum_{i=1}^{2}2u_iQ_i(x_0,x_1,x_2)+H(x_0,x_1,x_2),  \end{equation*}
with $u_1,u_2,x_0,x_1,x_2$ coordinates for $\Ps4$. The cubic threefold $X(Q,\theta)$ is smooth and contains the line $\bar\ell=\{x_0=x_1=x_2=0\}$.

Conversely, fix a pair $(X,\ell)$ with $X$ a smooth cubic threefold and $\ell$ a line on it.
Let $\pi_\ell$ be the projection of $X$ with center $\ell$ and $\tilde{\pi}_\ell: \tilde{X}\rightarrow \Ps2$ its resolution. The fiber over a point $p\in \Ps2$ is a conic $C_p$  coplanar with $\ell$. Let $Q(X,\ell)\subset \Ps2$ be the discriminant of this fibration, i.e., the  locus in $\Ps2$ parameterizing reducible conics. The curve $Q(X,\ell)$ is a smooth plane quintic for general $\ell$ in $X$ (cf. \cite[Lecture 2]{Tyu}).
Consider the curve 
\[ \tilde{Q}=\{(\ell',p) \in  G(2,5) \times Q(X,\ell) \mid \ell'\in \pi_\ell^{-1}(p), \ell\neq\ell' \}\] 
then $\tilde Q$ is an unramified double cover of $Q$. Let $\eta_\ell$ be the line bundle associated to this double cover, then $\eta_{\ell}\otimes \str_{Q}(1)$ is an odd theta characteristic $\theta_{\ell}$ on $Q(X,\ell)$ (cf. \cite{CMF}). 
We call the couple $(Q(X,\ell),\theta_{\ell})$ the {\em  discriminant curve} of the pair $(X,\ell)$.

With notation as above, the following holds:
\begin{Prop}\label{threefoldcons}
The spin curve $(Q,\theta)$ is the discriminant curve of the pair $(X(Q,\theta),\bar\ell)$.  \end{Prop}
In fact, it can be proven that the previous constructions give a birational correspondence (up to isomorphism on both sides): 
\[ \left\{ \begin{array}{c} \mbox{Smooth quintic plane curves with} \\ \mbox{a marked odd theta 
characteristic $\theta$} \end{array} \right\} \bircor \left\{\begin{array}{c} \mbox{Smooth cubic threefolds} \\ \mbox{with a marked line} \end{array} \right\}. \]

\begin{Rem}This result can be generalized to the case of a (reduced) nodal plane quintic plus an odd theta characteristic on its normalization (cf. \cite{CMF}). 

\end{Rem}
Consider a couple $(Q,\theta)$ and the corresponding matrix $M$ as above. The conic $T=\{L_{11}L_{22}-L_{12}L_{21}=0\}$ has an interesting geometric meaning:
\begin{Lem}\label{con}
The conic $T$ parameterizes plane sections $P=\bar\ell\cup C_p$ of $X$ with $C_p$ tangent to $\bar\ell$. In fact, it is the branch curve of $\tilde\pi_\ell|_E$ where $E$ is the exceptional divisor of the blow-up $\tilde X\ra X$. 

It is totally tangent to $Q$ and the intersection divisor $Q.T$ is the odd theta characteristic $\theta$ on $Q$.
\end{Lem}
\begin{proof}
The first two assertions follows from an easy calculation. The third one is proven in \cite[Proposition 4.2]{CMF}.
\end{proof}

Similarly, in the case of smooth quartics, we can associate to a couple $(C,\theta)$ a marked cubic surface $(X,p)$. 
\begin{Prop}\label{Cubicon} There is a one-to-one correspondence between
\[ \left\{ \begin{array}{c} \mbox{Quartic plane curves} \\ \mbox{with at most isolated singularities} \\ \mbox{and a marked bitangent} \end{array} \right\} \longleftrightarrow \left\{\begin{array}{c} \mbox{Irreducible cubic surfaces} \\ \mbox{with at most isolated singularities} \\ \mbox{and a marked smooth point} \end{array} \right\}  \]
up to isomorphisms on both sides.
\end{Prop}

\begin{proof} Suppose $C$ is a quartic plane curve and $b$ a bitangent of $C$. Let $Y'$ be the double cover of $\Ps2$ ramified along $C$, let $E$ be one of the pre-images of $b$ in $Y'$. Let $S'$ be the surface obtained by contracting $E$, then $S'$ is a cubic surface. The image of $E$ is the marked point $p$. (Contracting the other exceptional divisor gives an isomorphic cubic surface.)

Suppose that $S'$ is a cubic surface and $p\in S'$ a point. Let $\pi$ be the projection with center $p$. Resolving $\pi$ gives a morphism $Y' \ra \Ps2$ ramified along a quartic curve with at most isolated singularities. The marked bitangent is the image of the exceptional divisor of $Y'\ra S'$.

One can easily show that both constructions are each-other's inverse.
\end{proof}

\section{Del Pezzo surfaces and theta characteristics}\label{Comb}
In this section we generalize the notion of odd and even theta characteristic to quartic curves with at most $ADE$ singularities.

\begin{Not} For the rest of the section let $C$ be a quartic curve with at most $ADE$ singularities. Let $\pi':S'\ra \Ps2$ be the double cover of $\Ps2$ ramified along $C$. Let $S$ be the minimal desingularization of $S'$. Let $\pi: S \ra \Ps2$ be the composition of both maps. Let $\mathcal{N}\subset \Pic(S)$ be the subgroup generated by smooth rational curves with self-intersection $-2$. Let $\mathcal{L}_S$ and $\mathcal{BM}_S$ be the scheme of lines on $S$ and the scheme of blow-down models on $S$ (see Definition~\ref{Defsch}).
\end{Not}

On a smooth plane curve of degree 4 one can describe both the odd and the even theta characteristics in terms of bitangents. We prove that generalized bitangents on a stable quartic curve $C$  correspond in a natural way to odd spin structures on $C$, and that generalized Aronhold sets on $C$ correspond to even spin structures on $C$.

\begin{Def} 
A {\em generalized bitangent} on $C$ is a line $\ell \subset \Ps2$ such that either the intersection multiplicity $(C.\ell)_p$ is even at every point $p\in C\cap \ell$, or $\ell$ is a component of $C$.\end{Def}

\begin{Prop}\label{preimage} Suppose $E$ is a line on $S$  then $\pi(E)$ is a generalized bitangent  and the strict transform of ${\pi'}^{-1}(b)$ on $S$ consists of one or two lines. The first case occurs if and only if  $b$ is a component of $C$ or $b$ connects two double points of $C$.
\end{Prop}

\begin{proof} Cf. the proof of \cite[Proposition IX.1]{Dol}.\end{proof}

\begin{Def} Define the \emph{scheme of generalized bitangents} $\mathcal B(C)$ as $\mathcal{L}_S$ modulo the action of the Geiser involution.
\end{Def}

\begin{Rem} From Proposition~\ref{preimage} it follows that points of $\mathcal B(C)$ correspond to generalized bitangents. The reason to define $\mathcal B(C)$ in this way is to obtain the right multiplicity.

Another way of defining the multiplicity of a generalized bitangent would be to define an unramified degree 28 covering of the set of smooth quartic curves in the moduli space of GIT semi-stable quartic curves, and compactify this. This is the approach followed in \cite[Proposition 2.3.1]{CapSerbi}. It seems that taking the most natural compactification gives the same multiplicity as defined above. This strategy has a disadvantage: it does not define a multiplicity for bitangents on quartic curves that are not GIT semi-stable.
\end{Rem}

\begin{Prop}\label{prpoddcor} Suppose $C$ is a stable nodal plane curve of degree 4. Then there is an isomorphism from $S_C^{-}$ to  $\mathcal B(C)$. \end{Prop}

\begin{proof} We want to define an isomorphism of schemes $A:S_C^{-}\rightarrow \mathcal B(C).$
Let $\mathcal C\rightarrow T$ be a sufficiently general $1-$parameter family of quartics whose general fiber is smooth and whose fiber over $o\in T$ is $C.$
Let $S_{\mathcal C}^{-}$ and $\mathcal B(\mathcal C)$ be the corresponding families
of odd spin curves and bitangents.
On the generic fiber the correspondence between odd theta characteristics and bitangents is an isomorphism of reduced zero-dimensional schemes.
 So we get an isomorphism $\psi$ from $S_{\mathcal C}^{-}$ to  $\mathcal B(\mathcal C)$, at least away from the central fiber. As $S_{\mathcal C}^{-}$ is a smooth curve (see \cite[Proposition 2.2.1]{CapSertheta}), we can extend $\psi$ to the special fiber. 
Since $\psi$ is generically an isomorphism, it suffices to show that  $\psi|_{S_C^{-}}$ induces a bijection of sets between $S_{C}^{-}$ and $\mathcal{B}(C)$. It is explained in \cite[Section 3]{CapSerbi} 
how to construct this bijection.\end{proof}

It remains to generalize the notion of even spin curve. First of all, we generalize the notion of  Aronhold set.
\begin{Def} Let $C$ be a quartic plane curve with at most $ADE$ singularities.
A set of three distinct generalized bitangents $A=\{\ell_1,\ell_2,\ell_3\}$ is called {\em asyzygetic} if one of the following occurs:
\begin{itemize}

\item The intersection point of $\ell_i$ and $\ell_j$ is a singular point of $C$, for some $i\neq j$.
\item The  points of contact of the $\ell_i$ with $C$ do not lie on a conic $K$ such that all the hyperflex lines of $C$, which are contained in $A$, are tangents of $K$.

\end{itemize}\end{Def}

We recall the definition of Aronhold set for a smooth quartic curve:
\begin{Def}
Suppose $C$ is smooth. A set $\{\ell_1,\dots, \ell_7\}$ of seven bitangents  is called an \emph{Aronhold set} if for all $I=\{i,j,k\} \subset \{1,2,\dots,7\}$ such that $\# I =3$ the triple $\{\ell_i,\ell_j,\ell_k\}$ is asyzygetic.
\end{Def}

For our aims it will not suffice to define an Aronhold set in terms of generalized bitangents. If $C$ is smooth then there is a 2:1 correspondence between Aronhold sets on $C$ and blow-down models of $S$. The data of generalized bitangents together with a multiplicity will not always determine uniquely a blow-down linear system up to the Geiser involution. See Example~\ref{exaunique}.

\begin{Def}  
The {\em ADE-zation}  of $C\subset \Ps2$ is a triple $(C',R,\psi)$ with $R$ a surface, $\psi: R \ra \Ps2$ a morphism and $C'=\psi^{-1}(C)$. The morphism $\psi$ is such that $\psi^{-1}(p)$ is a point if $p\in \Ps2\backslash C_{\sing}$  and  consists of a union of smooth rational curves $E_{i,p}$ if $p\in C_{\sing}$. 
For every $p\in C_{\sing}$ the intersection-numbers  $E_{i,p}.E_{j,p}$, are as on \cite[page 88]{BPV}. 
\end{Def}

The existence of the ADE-zation can be proven as follows. 
Let $\iota$ be the involution associated to the double cover $\pi'$. We can lift $\iota$ to $S$ in an unique way, call this involution $\sigma$. Then $R=S/\langle \sigma \rangle$ is a smooth surface, let $p: S\ra R$ be the induced morphism. Then the reduced scheme structure on $p_*\pi^*(C)$ gives the ADE-zation curve $C'$. The morphism $\psi$ is the unique morphism such that $\psi p = \pi$.

\begin{Def}\label{defgenAro}
Let $(C',R,\pi)$ be the ADE-zation of $C$.  An {\em Aronhold line} is a smooth rational curve $D$ on $R$ that is either  the strict transform of a generalized bitangent of $C$ (which we call {\em type 1}) or a component of $C'$ not contained in the strict transform of $C$ (which we call {\em type 2}).

An (unordered) 7-uple $A=D_1+\dots+ D_7$ of  Aronhold lines is called {\em a generalized  Aronhold set} if the following  conditions hold:
\begin{enumerate}
\item Suppose $D$ is an Aronhold line occurring more then once in $A$, then $D$ is of type 2, occurs exactly twice in $A$ and  either $D^2=-2$ or $D^2=-1$ and $D$ is tangent at the strict transform of $C$.
\item Every  connected component $\Gamma$ of the dual graph of the Aronhold lines of type 2 (without multiplicities) in $A$, is of type $A_k$, for some $k>0$. 
\item For every connected component $\Gamma$ of the dual graph of the Aronhold lines of type 2 in $A$,  we have an unique Aronhold line $D\in A$ of type 1, which intersects some of the components corresponding to the vertices of $\Gamma$. Moreover, $D$ intersects at most one component, and intersects it transversally. 
\item For every connected component $\Gamma$ of the dual graph of the Aronhold lines of type 2 in $A$, we have that the subgraph $\Gamma_2$  of Aronhold lines occuring twice in $A$, is connected, and either a $-1$-curve contained in $C$ intersects one of the components corresponding to $\Gamma$ or one of the Aronhold lines occuring twice intersects the strict transform of $C$.
\item Any three distinct Aronhold lines of type 1 in $S$ are asyzygetic.
\end{enumerate}
\end{Def}

\begin{Rem} Suppose $C$ is smooth, let $A$ be a generalized Aronhold set, then $A$ is a classical Aronhold set.\end{Rem}

\begin{Exa}\label{exaunique} Suppose $C$ has a singularity of type $A_5$ and one of type $A_2$. Let $D$ be the strict transform on $R$ (the ADE-zation surface) of a  bitangent passing through the singular points. Let $E_1,E_2, E_3$ be the Aronhold lines of type 2, at the singular point of type $A_5$, (where $E_3$ intersects the strict transform of $C$) and $E'_1$ is the Aronhold line of type 2 over the cusp.

Then both 
\[  D+2 E_1+2 E_2+E_3+E'_1 \} \mbox{ and } D+E_1+2 E_2+E_3+2E'_1\] are generalized Aronhold sets, hence $D$ and its multiplicity do not determine the generalized Aronhold set.

In a very similar way one can construct singular cubic surfaces $Y$ such that one limit position of a double-six (i.e. lines with a given multiplicity) corresponds to more than one point in $\mathcal{D}_{\tilde{Y}}$, where $\tilde{Y}$ is the minimal desingularization of $Y$.
\end{Exa}

\begin{Prop}\label{syzPrp} Let $A=\{b_1,\dots, b_n\}$ be a set of Aronhold lines of type 1. Then the following are equivalent:
\begin{enumerate}
\item Any triple of distinct Aronhold lines is asyzygetic.
\item We can choose lines $L_i, i=1,\dots, n$ on $S$, such that $L_i.L_j=0$ if $i\neq j$, and $\pi(L_i)=b_i$.
\end{enumerate}\end{Prop}

\begin{proof} Recall the following diagram, where $R$ is the ADE-zation surface:
\[ \begin{array}{ccc}
 S & \stackrel{p}{\ra} & R \\
\downarrow  &\searrow \pi & \downarrow \psi \\
  S' & \stackrel{\pi'}{\ra} & \Ps2 \end{array}.\]
We have that (2) is not satisfied if and only if $\pi|_{\pi^{-1}(b_i\cup b_j \cup b_k)}$ admits a section for some $i,j,k$ pairwise distinct. This is equivalent to ${\pi'}|_{{\pi'}^{-1}(b_i \cup b_j \cup b_k)}$ admits a section and none of the points in $b_i\cap  b_j$, $b_j\cap b_k$ and $b_j \cap b_k$ are in the singular locus of $C$. A reasoning as in \cite[Proposition III.1.7]{Ver} shows that we can find a section to ${\pi'}|_{{\pi'}^{-1}(b_i\cup b_j \cup b_k)}$, only when the points of contact of the $b_i$ with $C$ lie on a conic $K$ such that if $b_i$ is hyperflex line then it is tangent to $K$. From this we obtain the proposition.
\end{proof}

\begin{Prop}\label{PrpArocnd} Let $A=b_1+\dots+ b_7$ be an unordered 7-uple of Aronhold lines. Then the following are equivalent:
\begin{enumerate}
\item $A$ is a generalized Aronhold set.
\item There are 7 exceptional curves $E_1,\dots, E_7$ in $S$, such that we can blow down subsequently $E_1,\dots, E_7$, if $D\in A$ then at least one of the $E_i$ is contained in $\pi^{-1}(D)$ and if $D$ occurs twice in $A$ then there are $i,j, i \neq j$, such that $\pi^{-1}(D)=E_i\cup E_j$.
\end{enumerate}\end{Prop}

\begin{proof} Suppose $A$ is a generalized Aronhold set. Then the pre-image of an Aronhold line $b_i$ of type 1 consists of one irreducible curve $D_i$ or two irreducible curves  $D_{i,1}$ and $D_{i,2}$, with self-intersection $-1$. The fifth condition in Definition~\ref{defgenAro} and  Proposition~\ref{syzPrp} imply that we can label the lines over $\pi^{-1}(b_i)$ as $D_{i,1}$ and $D_{i,2}$ such that $D_{i,k}$ and $D_{j,k}$ do not intersect whenever $i\neq j$ (if $\pi^{-1}(b_i)$ consists of one line then we take $D_{i,1}=D_{i,2}=E_i$ and the intersection property holds also). If $b_i$ is of type 1, we take $E_i:=D_{i,1}$.

Suppose $b_i$ is an Aronhold line of type 2 occuring twice in $A$. Then choose $E_i, E_j$ such that $\pi^{-1}(b_i)=E_i\cup E_j$.

Consider now the $b_i$ which are Aronhold lines of type 2 occuring only once in $A$. We have a unique way of choosing exceptional curves $E_i$ over $b_i$, such that each connected component of the dual graph of the $-2$ curves in $\{E_1,\dots, E_7\}$ intersects a unique line contained in $\{E_1,\dots, E_7\}$. Hence we can reorder $E_1,\dots, E_7$ to obtain a blow-down morphism to $\Ps2$.

Conversely, if (2) holds, then the first four conditions of Definition~\ref{defgenAro} are satisfied, (see for example \cite[page 88]{BPV}), the fifth condition follows from Proposition~\ref{syzPrp}.\end{proof}

\begin{Cor}\label{Corlinaro} Let $C$ be a quartic curve with at most $ADE$ singularities. Let $S$ be the desingularization of the double cover of $\Ps2$ ramified along $C$. Then there is a 2:1 correspondence between blow-down models $S \ra \Ps2$ and generalized Aronhold sets.\end{Cor}

\begin{Rem} In \cite[Section 7]{GrHa} there is a 1:1 correspondence between blow-down models $S\ra\Ps2$ and Aronhold sets, but in that paper the blow-down linear systems are considered up to the Geiser involution.\end{Rem}

\begin{Def} Define the {\em scheme of Aronhold sets} $\mathcal A(C)$ as the quotient of $\mathcal{BM}_S$ by the Geiser-involution.

Let $b_0$ be a generalized bitangent. Let $E$ be a smooth irreducible curve on $S$ such that $\pi(E)=b_0$. Let $\mathcal{BM}_{S,E}$ be the scheme of blown-down models $|L|: S \ra \Ps2$ such that $L.E=0$. The multiplicity of $|L|$ equals
\[\frac{\# \{ \sigma \in \Aut(NS(S)) : \sigma(E)=E, \sigma(L)\equiv L \bmod \mathcal{N}\} }{\# \{ \sigma \in \Aut(NS(S)) : \sigma(E)=E, \sigma(L)=L \}}. \]
 
The {\em scheme of Aronhold sets containing $b_0$} is $\mathcal{BM}_{S,E}$ modulo the Geiser involution and is denoted by $\mathcal A(C,b_0)$. On $\mathcal A(C,b_0)$ acts the involution induced by the involution of $\Pic(S)$ fixing the irreducible components over $b_0$. The quotient $\mathcal{T}(C,b_0)$ by this involution is the {\em scheme of generalized even theta characteristics} on $(C,b_0)$.
\end{Def}

\begin{Prop}\label{doublesix} 
Let $b_0$ be a generalized bitangent of $C$ not passing through any of the singular points of $C$.  Let $(X,p)$ be the marked cubic surface corresponding to $C$ from Proposition~\ref{Cubicon}. 
Then there is an isomorphism from $\mathcal A(C,b_0)$ to $\mathcal{BM}_X$. 

Let $X'$ be the nodal Del Pezzo surface corresponding to $X$. Then there is an isomorphism between $\mathcal{D}_{X'}$ and $\mathcal T(C,b_0)$. 
If $C$ is stable the two schemes are both isomorphic to $S_C^{+}$.
\end{Prop}

\begin{proof} Our assumption on $b_0$ implies that $p$ is a smooth point on $X$, not on any of the lines of $X$. 

From Proposition~\ref{PrpArocnd} it follows that there is a $2:1$ correspondence between blow-down linear systems on $S$ (which is the desingularization of $X$ blown up in $p$) and generalized Aronhold sets on $C$. From this we obtain a $2:1$ correspondence between blow-down linear systems on $X'$ and Aronhold sets containing $b_0$. If $X$ were smooth then there would be a $1:1$ correspondence between blow-down models on $X$ and Aronhold sets containing $b_0$ (see Remark~\ref{ds}). By definition of $\mathcal A(C,b_0)$ this correspondence extends to the singular case.

Suppose now that $C$ is a stable curve. Let $\mathcal{C}\ra T$ be a (sufficiently general) 1-parameter family of plane quartics, such that the generic fiber is smooth and the fiber over $o\in T$ is $C$. Let $\mathcal{S}_\mathcal{C}^{+}$ be the scheme of even spin structures on $\mathcal{C}$, let $\mathcal{T}({\mathcal C},\mathcal{B}_0)$ be the corresponding scheme of Aronhold sets containing $\mathcal{B}_0$ modulo the involution. There is an isomorphism of schemes away from the central fiber. Since $\mathcal{S}_\mathcal{C}^{+}$ is a smooth curve, it suffices to give a bijection of points on the central fiber, to obtain an isomorphism $\mathcal{S}_\mathcal{C}^{+}\cong \mathcal{T}(\mathcal{C},\mathcal{B}_0)$.

Given a generalized Aronhold set $B$, we have $B=\{E_1,\dots, E_n,F_1,\dots, F_{7-n}\}$, where the $E_i$ are Aronhold lines of type 1 and the $F_j$ are (distinct) Aronhold lines of type 2. We denote with $F(E_j)$ the (unique) Aronhold line in $B$ of type 1 intersecting $E_j$. Let $\theta$ be the line bundle on $C'$, the ADE-zation curve of $C$, such that 
\[ \theta= \str_{C'}\left(\frac{1}{2} \left(\sum_i F_i+\sum_j F(E_j)+2 E_j-6\bar{\ell}\right) \cap C'\right),\]
with $\bar{\ell}$ the strict transform of a general line in $\Ps2$. Then on every exceptional component of $C$ the degree of $\theta$ is 0 or 1. Let $Y$ be the curve obtained by contracting all the components of $C'$ such that $\theta$ restricted to that component has degree 0. Then the push forward of $\theta$ gives a line bundle on $Y$. Since $(Y,L)$  is the limit of family of even spin structures, we obtain that $(Y,L)$ is  even.

Given a spin curve $(Y,L)$, take a family of smooth plane quartics with even theta characteristics with central fiber $(Y,L)$. Associate to this an Aronhold set containing $b_0$ (this is possible away from the central fiber). Let $D_1, \dots, D_7$ the limit position of the bitangents in the Aronhold set. These lines are of the form $E_i$ (of type 1) or a sum of Aronhold lines of type 1 and type 2, with only one of type 1. The set $A$ consisting of the irreducible components of the $D_i$ has 7 elements. Since being asyzygetic is an open condition, all triples of Aronhold lines of type 1 in $A$ are asyzygetic.

Since $b_0$ does not pass through any singular points of $C$ we have a canonical isomorphism $\mathcal{BM}_{X'} \ra \mathcal{BM}_{S,E}$, where $E$ is one of the components of the pre-image of $b_0$. This induces an isomorphism $\mathcal{T}(C,b_0)$ and $\mathcal{D}_{X'}$.
\end{proof}

\begin{Rem} If we drop the assumption that $b_0$ does not pass through any singular point, then there exists at least two distinct double-six on $X$ giving rise to the same point in $\mathcal{T}(C,b_0)$.
\end{Rem}

\begin{Rem}\label{mincur} We would like to use points in either $\mathcal A(C)$ or $\mathcal A(C,b_0)$ to define a line bundle on the ADE-zation curve $C'$.

Fix a point $\theta$ in $\mathcal A(C)$. Fix a blow-down linear system $L\in \mathcal{BM}_S$ over $\theta$. 
We can find seven distinct, effective and reduced divisors $D_i$ on $S$ such that
\begin{itemize}
\item The support of $D_i$ contains exactly one line $E(i)$, and $E(i)$ is contracted by $L$.
\item All the irreducible components  of $D_i-E(i)$ are smooth rational curves $F_j$, with $F_j^2=-2$ and $F_j.L=0$.
\item The dual graph of the support of $D_i-E(i)$ is connected.
\item $D_i.D_j=-\delta_{i,j}$.
\end{itemize}
The most natural line bundle to associate to $\theta$ is the intersection divisor of $C'$ with the sum of the push-forwards of the $D_i$ minus 3 times a general line  and this divided by two  (as in the proof of Proposition~\ref{doublesix}). If we choose a different blow-down linear system, corresponding to the same Aronhold set,  we obtain 7 divisors on $S$ which differ from the $D_i$ by the Geiser involution.\end{Rem}

\begin{Exa}\label{node} Consider a family of Del Pezzo surfaces of degree 2. Suppose that the generic fiber is a smooth Del Pezzo surface and the special fiber $S$ is a nodal Del Pezzo surface with a unique $-2$-curve $F$. 

The linear system $|-K_S|$ is the blow-down of $S$ onto a normal surface $S'$ with a node combined with a degree 2 morphism $\pi: S \ra \Ps2$. The morphism $\pi$ is branched along a quartic curve $C(S)$.
One way of obtaining such a family is to take $7$ points in general position $P_i$ in $\Ps2$ and then move $P_3$ to a point on the line connecting $P_1$ and $P_2$ with $P_3\neq P_1,P_2$. We consider now this case.
The Picard group of $S$  has still rank 8 and is generated by $L$ and the $E_i$. In Notation~\ref{mark} we gave 56 divisors expressed in terms of $L$ and the $E_i$. Of these divisors only 44 define irreducible divisors on $S$. The reducible ones are the $L_{i,j}$ and $D_i$ for $1\leq i<j \leq 3$ and $C_{i,j}$ for $4\leq i<j\leq 7$. They satisfy the following relations in $\Pic(S)$:
\begin{eqnarray}\label{picrela} L_{i,j}=F+E_k,  D_k=F+C_{i,j} &\mbox{ for }& \{i,j,k\}=\{1,2,3\}, \\ \label{picrelb} C_{i,j}=F+L_{k,m}&\mbox{ for } & \{i,j,k,m\}=\{4,5,6,7\}.\end{eqnarray}
In total we have 56 curves on $S$ with self-intersection $-1$. If we contract $F$ then we get 12 pairs of divisors that coincide.

The curve $C(S)$ is a curve with one node $P$. If $D$ is an effective divisor on $S$ with $D^2=-1$ then $D.S=0$ if and only  if $b_D$ does not pass through the node $P$. In total there are 16 bitangents not passing through the node, and 6 bitangents passing through the node. The latter 6 bitangents correspond to theta characteristics with multiplicity 2. 

For the even theta characteristics, we study the Aronhold sets, or, equivalently  the  linear systems on $S$ giving blow-downs $S\ra \Ps2$, modulo the Geiser involution.
We characterize the degeneration of the linear system $| D| $ by the intersection number $F.D$:
\begin{enumerate}
\item Suppose $F.D=2$, then $| D| :S\ra \Ps2$ is the blow-up of seven points $P_1,\dots, P_7$ such that $P_1,\dots, P_6$ are on a conic,  no other set of six points lie on a conic, no three points are on a line and no two points coincide.
\item Suppose $F.D=1$, then $| D| :S\ra \Ps2$ is the blow-up of seven points $P_1,\dots, P_7$ such that $P_1,P_2$ and $P_3$ are collinear,  no other set of three points lie on a line, no six points are on a conic and no two points coincide.
\item Suppose $F.D=0$, then $| D| :S\ra \Ps2$ is the blow-up of seven points $P_1,\dots, P_7$ such that $P_1$ and $P_2$ are infinitesimal close,  no other two points are infinitesimal close, no three points lie on a line and no six points are on a conic.
\item Suppose $F.D=-1$, then $F$ is in the fixed part of $| D |$, and the movable part of $| D |$ is a linear system $ | \tilde{D}|$ with $\tilde{D}.F=1$.
\item Suppose $F.D=-2$, then $F$ is in the fixed part of $| D |$. The linear system $\tilde{D}=D-2F$ is a blow-down linear system.

\end{enumerate}
The frequencies of intersection numbers are given below:
\[ \begin{array}{l|lllll}
\mbox{Linear system $| D| \; \backslash$ } D.F  & 2 & 1 & 0 &-1& -2 \\
\hline
L &  0&1&0&0 &0 \\
2L-E_m-E_n-E_p& 4 & 18 & 12 &1&0\\
3L-\sum_{t=1}^7 E_t +E_i+E_j-E_k &12 & 39 &36 & 18 & 0\\
4L-\sum_{t=1}^7 E_t +E_i -E_m-E_n-E_p &12 & 40 &48 & 36 & 4 \\
5L-2 \sum_{t=1}^7 E_t + 2E_i &0&3&0&4&0 \\
8L-3\sum_{t=1}^7 E_t & 0&0&0&1 &0\\
7L-2\sum_{t=1}^7 E_t-E_m-E_n-E_p &0&1&12&18&4\\
6L-2\sum_{t=1}^7 E_t -E_i-E_j+E_k & 0&18&36&39&12\\
5L -\sum_{t=1}^7 E_t-2E_i-E_j-E_k-E_l&4&36&48&40 & 12\\
4L-\sum_{t=1}^7 E_t -2 E_i & 0&4&0&3&0\\
\hline
\mbox{total} & 32 & 160 & 96 &160 & 32 \\
\end{array}.\]
To each even theta characteristic we can associate 16 Aronhold sets. 
The linear systems with $D.F=\pm 1$ identify all Aronhold sets on $C(S)$ for which the associated even theta characteristic has multiplicity 2. 
So $\frac{160}{16}=10$ even theta characteristics have multiplicity 2. 

It turns out that for every even theta characteristics with multiplicity 1 there are exactly 2 linear systems $|D|$ with $D.F=2$, there are 2 linear systems with $D.F=-2$ and 12 linear systems with $D.F=0$. In total there are $\frac{256}{16}=16$ even theta characteristics with multiplicity 1.
The symmetry in the above table follows from the fact that the Geiser involution sends $F$ to $-F$.
\end{Exa}
\begin{Exa}\label{cusp}
Consider a family of Del Pezzo surfaces of degree 3. Suppose that the generic fiber is a smooth Del Pezzo surface and the special fiber $S$ is a degenerate Del Pezzo surface with two -2-curves $F_1, F_2$ intersecting transversally. 

The linear system $|-K_S|$ is the blow-down of $S$ onto a cubic surface with an $A_2$ singularity combined with a degree 2 morphism $\pi$ branched along a quartic curve $C(S)$ with a cusp.
One way of obtaining such a family is to take $6$ points in general position $P_i$ in $\Ps2$ and then move $P_1$ infinitely close to $P_2$ and $P_3$ infinitely close $P_2$. 
Define $E_2=F_2+E_1$ and $E_1=F_1+F_2+E_3$.

Then $E_3$ up to $E_6$ still exists. Define $E_2=F_2+E_1$ and $E_1=F_1+F_2+E_3$. Then $L_{1,k}$, for $k=2$ and $k\geq 4$, the line $L_{i,j}$ for $4\leq i \leq j \leq 7$ and $C_k$ for $k\geq 3$ are irreducible. The other $-1$-curves satisfy the following relations 
\begin{eqnarray}\label{picrela2} L_{2,k}=L_{1,k}+F_1,\; L_{3,k}=L_{1,k}+F_1+F_2, & \mbox{ for } k\geq 4 \\ 
C_2=C_3+F_2, \;  C_1=C_3+F_1+F_2 \\
L_{1,3} = L_{1,2} +F_2, \; L_{2,3}=L_{1,2}+F_1+F_2 \end{eqnarray}
In total there are six $-1$ curves with multiplicity 3, and nine $-1$ curves with multiplicity 1. From this we can deduce that a quartic curve with a cusp has 10 proper bitangents and six bitangents passing through the cusp.

There are 12 linear systems with multiplicity 1, there are 18 linear systems with multiplicity 3, and one linear system $L$ with multiplicity 6. The linear system $L$ gives rise to the double six $(L,L)$. This double-six has multiplicity 3. Hence on a quartic curve with a cusp there are 6 even theta characteristics of multiplicity 1, and 10 of multiplicity 3.
\end{Exa}

\begin{Exa}\label{E_7} Suppose $S$ is a Del Pezzo surface of degree 2 containing 7 curves of self-intersection $-2$ such that they intersect as $E_7$. 
Then $X$ contains a unique line, so $C(S)$ has only one  generalized bitangent $b$, and $b$ is the unique line intersecting $C(S)$ with multiplicity 4. There is an unique generalized Aronhold set $A$ and a unique blow-down linear system.

The unique Aronhold set $A$  has multiplicity 36 in $\mathcal A(C,b)$, and multiplicity 288 in $\mathcal A(C)$. The fact that the multiplicity of $A$ is different in both schemes comes from the fact that $b$ has a multiplicity in $\mathcal B(C)$.
\end{Exa}

\section{Clemens' model}\label{oldmodel}
An essential tool in the study of the geometry of the cubic threefold  $X\subset \Ps4$ is given by its \emph{intermediate Jacobian} $J(X)$, defined as
\[ J(X)=({H^{3,0} \oplus H^{2,1}})^{*}/H_3(X,\Z).\]
Since $H^{3,0}(X)=0$ we have that $J(X)$ is a principal polarized abelian variety.
We denote with $\Theta\subset J(X)$ the theta divisor of $J(X)$.
Mumford has established a connection between the intermediate Jacobian of cubic threefolds and Prym varieties. Given a smooth cubic threefold $X$ and a line $\ell \subset X$, there is an isomorphism of principally polarized abelian varieties:
\[(J(X),\Theta)\cong (\Pry(C,\eta),\Xi),\]
where $(\Pry(C,\eta),\Xi)$ is the Prym variety associated to the double cover $\eta:\tilde C\rightarrow C$ coming from the conic bundle induced by $(X,\ell)$ (see Section~\ref{beau}). 
Analogous to the case of Jacobians of smooth projective curves, Clemens and Griffiths defined an
 \emph{Abel-Jacobi map} (see \cite[Section 4]{CG}):
\[ A:A_1(X)\rightarrow J(X),\]
where $A_1(X)$ is the group of algebraic 1-cycles on $X$ homologically equivalent to zero modulo the group generated by cycles rationally equivalent to zero. For a family $\{\Gamma_b\}_{b\in B}$ of algebraic 1-cycles with a base point $b_0$, denote by $A_B$ the Abel-Jacobi map base changed by the map $b\mapsto \Gamma_b-\Gamma_{b_0}$.
In this section we consider the cases where $B$ is either the Fano scheme $F(X)$ of lines on $X$ or the scheme $T$ of rational cubics on $X$. The maps $A_{F(X)}$ and $A_T$ have interesting geometric interpretations and both lead to parameterizations of the theta divisor $\Theta\subset J(X)$.

A complete description of the Abel-Jacobi map $A_F$ on the Fano scheme is due to Clemens and Griffiths (cf. \cite{CG}). In the following, let $\psi:(F\times F)-\Delta \rightarrow {\Ps4}^*$ be the map:
\[ (\ell,\ell')\mapsto \left\{  \begin{array}{c}
\spa(\ell,\ell')\mbox{ if } \ell\cap \ell'=\emptyset\\
T_{x}X\mbox{ if } \{ x\} = \ell\cap\ell' \end{array}\right.\]
Let $\G: \Theta \ra{ \Ps4}^*$ be the Gauss map, via identification of $\Ps{} (H^{1,2}(X))$ with ${\Ps4}^*$ through Griffiths residue calculus (see \cite{Gr}).
Let $\phi: F \times F \rightarrow J(X)$ be the difference of the Abel-Jacobi map, i.e. $\phi(\ell,\ell')=A_F(\ell)-A_F(\ell')$.
\begin{Thm}[Clemens-Griffiths {\cite{CG}}]\label{gauss} We have the following
\begin{itemize}
\item The image of $\phi$ is contained in $\Theta$ and $\phi$ has degree 6 onto its image. 
\item The general fiber of $\phi$ is given by a double six $\{(l_1,m_1),\dots,(l_6,m_6)\}$ ($l_i$ pairwise skew, $m_i$ pairwise skew, $l_i$ and $m_j$ skew if and only if $i=j$) in a smooth hyperplane section of X.
\item On $F\times F-\Delta$ we have $\G\phi=\psi$.
\item The branch locus of $\G$ equals $X'$, the dual variety of $X$.
\end{itemize}
\end{Thm}
A different model of the theta divisor using smooth cubic curves on $X$ has been given by Clemens in \cite[Section 4]{C}. Cubic curves in the family $T$ are of two types: those contained in a plane sections and smooth twisted cubics contained in an hyperplane. 
It can be easily proven that if $X\cap H$ is a smooth hyperplane section, then the smooth twisted cubics in $H$ give all the blow-down linear systems on the cubic surface $X\cap H$. 

\begin{Thm}[Clemens {\cite{C}}]\label{model} We have the following:
\begin{itemize}
\item The image of the Abel-Jacobi map $A_T:T\rightarrow J(X)$ is the theta divisor $\Theta\subset J(X)$.
\item The general fiber of $A_T$ consists of a linear system $L_H$ of twisted cubics in a single hyperplane section $X_H=H\cap X$. In particular, this holds for hyperplane sections of $X$ with $A_1$-singularities.
\item The set of rational plane cubics is mapped by $A_T$ to the only singular point of $\Theta$.

\end{itemize}
\end{Thm}
Let $\mathcal U$ be the set of hyperplane sections having at most rational double points. Let
\[\mathcal{C}=\{(H,L_H):H\in \mathcal{U}, L_H \in \mathcal{BM}_{X_H} \}, \]
with, in the case that $X_H$ is singular, we define $\mathcal{BM}_{X_H}$ as $\mathcal{BM}_{\tilde{X}_H}$ with $\tilde{X}_H$ the desingularization of $X_H$.
\begin{Cor}
The Abel-Jacobi map $A_T$ induces a birational correspondence between $\mathcal{C}$ and $\Theta$.
\end{Cor}
Combining the previous results, we obtain a simple description of the Gauss map on the theta divisor:
\begin{Prop}\label{gaussbir}
There exists a birational map $\mathcal A:\mathcal {C}\dashrightarrow \Theta$ such that $\mathcal G \mathcal A=\pi$ where $\pi:\mathcal {C}\rightarrow {\Ps4}^*$ is the projection on the first factor. Moreover the map $\mathcal A$ induces a birational isomorphism $\tilde{\mathcal A}$ between 
\[ \tilde{{C}}=\{(H,D_H):H\in \mathcal{U},\ D_H \in \mathcal D_{X_H}\} \dashrightarrow \tilde \Theta= \Theta /\langle -1\rangle.\]
\end{Prop}
In the following we call $\tilde{C}$ {\em Clemens' model} for the theta divisor $\Theta$. Essentially, Clemens' model is the blow-up of $\Theta$ at its triple point.

\section{The new model}\label{newmodel}
Let $X$ be a smooth cubic threefold in $\Ps4$. Assume that all isolated singularities on  hyperplane sections $X\cap H$  are of type $ADE$. Fix a line $\ell$ on $X$. Let $(Q,\theta)$ be the associated pair from Proposition~\ref{threefoldcons}, where $Q$ is a quintic plane curve and $\theta$ an odd theta characteristic on $Q$. Assume that $Q$ is smooth and  $h^0(\theta)=1$ (the last assertion is equivalent to $\theta \not \cong \str_{\Ps2}(1)|_Q$). For  a general smooth cubic threefold $X$ the above assumptions are satisfied. We identify divisors and line bundles, whenever no confusion arises.

Let $V=|  \theta(1) | $ and $d\in V$. Since 
 \[2d\in | \theta(1)^{\otimes 2}| =| \str_{Q}(4)| \]
and $d$ is effective, we have that $2d$ is cut out by a unique quartic $C_d\subset\Ps2$. By varying $d$ we obtain a four dimensional irreducible family of quartics totally tangent to $Q$. Note that $d$ is also a divisor on $C_d$. If $C_d$ is smooth, then the divisor $\theta_d:=3K_{C_d}-d$ gives an (odd) theta characteristic on $C_d$, call the corresponding bitangent $b_d$. 
 The idea is to take the family of quartics in $V$ (possibly singular) with all possible even theta characteristics on them.

\begin{Not} 
Let $U\subset V$ be the set of $d\in V$ such that $C_d$ is non-reduced.
Let $B^0=\{(C_d,b_d)\mid C_d \mbox{ smooth}\}$ and $B$ its closure in $(V-U)\times {\Ps2}^*$. We define:
\[V_{Q,\theta}:=\{ (C_d,b_d,A) | (C_d,b_d)\in B , \; A\in \mathcal{T}(C_d,b_d)\}.\]
\end{Not}
  
Consider the natural projection
\[\mathcal F:V_{Q,\theta}\rightarrow V \ \  (C_d,b_d,A)\mapsto d.\]
Since a smooth quartic possesses 36 even theta-characteristics, $\mathcal{F}$ is generically $36$-to-$1$ and hence it is dominant. We show in Proposition~\ref{strangecon} that the elements of $U$ are conics with a double line.
In the sequel, we show that $V_{Q,\theta}$ is a birational model for the quotient
$\tilde\Theta=\Theta/\langle -1 \rangle$ and the map $\mathcal F$ is essentially the
Gauss map $\tilde{\G}:\Theta/\langle -1 \rangle\rightarrow
{\Ps4}^*$. 
\begin{Def}\label{Prymcan} 
Let $\eta:=\theta(-1)\in \Jac(Q)[2]$. The morphism defined by $V=| K_Q\otimes \eta |$ is called  \emph{Prym canonical map} associated to $\eta$:
\[\phi:Q\rightarrow Q'\subset |K_Q\otimes \eta| ^{*}.\]
We can define an isomorphism $\Lambda^* : \mid K_Q\otimes\eta\mid^* \ra \Ps4$ such that $\{ \Lambda^*(\varphi(x)) \}$ is the singular point of $\pi_\ell^{-1}(x)$ (cf. \cite{Beau81}).

 If $d\in\mid K_Q\otimes\eta\mid$ then we denote by $H_d$ the corresponding hyperplane in $\Ps4$.
We call a point in the intersection of $X \cap Q'$ a {\em vertex}. 
We call $\Lambda: V \ra {\Ps4}^*$ the {\em dual Prym canonical map}. Denote $X_d := X\cap \Lambda(d)$ for $d\in V$. 
\end{Def}

\begin{Lem}\label{qua} Let $d\in V$ and let $\ell$ be a line not contained in $H_d$. Then the quartic $C_d$ is the ramification curve of the projection $\pi_d$ of $X_d$ from $p_d=H_d\cap \ell$.
\end{Lem}
\begin{proof}
Let $B_d$ be the ramification curve of the projection $\pi_d$. Notice that $\pi_d=\pi_{\ell}|_{X_d}$. A point $x\in \Ps2$ lies on the intersection $B_d\cap Q$ if and only if the hyperplane $H_d$ cuts the plane spanned by $\ell_x$ and ${\ell'}_x$ in a line containing the vertex $v=\ell_x\cap \ell'_x$, where $\ell_x \cup {\ell'}_x=\pi_{\ell}^{-1}(x)$.

On the other hand, the support of $d$ is given by ${\Lambda^*}^{-1}(Q'\cap H_d)$, that is, by planes with vertexes on $H_d$. Then the quartic $B_d$ cuts the quintic $Q$ in $d.$ 
Since $C_d$ and $B_d$ cut out the same divisor $d$ on $Q$ it follows that $B_d=C_d$.
\end{proof}

\begin{Prop}\label{Phi}
There is a birational map $\mathcal B: V_{Q,\theta}\rightarrow \tilde\Theta$.
\end{Prop}

\begin{proof}It suffices to give a birational correspondence $\Phi$ between $V_{Q,\theta}$ and Clemens' model $\tilde{\mathcal{C}}$. 

The correspondence between $V$ and ${\Ps4}^*$ is given by $\Lambda$.
Suppose $H_d$ is a hyperplane in ${\Ps4}$ such that $\ell$ is not contained in $H$. Then Lemma~\ref{qua} assures that the corresponding plane quartic $C_d$ is the ramification curve of the resolution of $\pi_\ell|_{H_d}$. Moreover, if $H$ is a generic hyperplane then Proposition~\ref{doublesix} implies that the double-six on $X_d$ correspond to the even theta-characteristics on $C_d$. This defines the birational correspondence between $V_{Q,\theta}$ and $\tilde{C}$.
\end{proof}

\begin{Cor}
Let $\tilde{\mathcal{G}}:\tilde{\Theta}\rightarrow {\Ps4}^*$ be the Gauss map. Then $\tilde{\mathcal{G}} \mathcal B=\mathcal F.$
\end{Cor}

In the following we describe  the ramification divisor of the  model $\mathcal{F}:V_{Q,\theta}\rightarrow V$ for the Gauss map $\tilde{\mathcal{G}}$. 
Moreover, we compare it with Clemens' model $\tilde\pi:\tilde C \rightarrow {\Ps4}^*$.
We have seen (Proposition~\ref{gauss}) that the branch locus of $\mathcal G$ is given by the dual variety of $X$. The branch locus of $\mathcal F$ has additional components.  

\begin{Not} Denote
\[ \overline{C_{\ell}}:=\{(H,D_H) \in\tilde{C} \mid \ell\subset H\} \mbox{ and } \mathcal{H}:=\{(H,D_H) \in \tilde{C} \mid H\in X'\} \]
It follows from the proof of Proposition~\ref{Phi}  that we obtain a birational  map
\[ \Psi: \tilde{C} \dashrightarrow V_{Q,\theta}.\]
\end{Not}

The map $\Psi$ is  defined on an open dense subset $U'$ of $\mathcal{H}$. For generic $(H,D_H)\in U'$ the curve $C_H$ is a quartic with one node $q$, and $q$ does not lie on the quintic $Q$. 

Let $\mathcal H_1$ be the closure in $V_{Q,\theta}$ of the locus corresponding to quartics with a node not on $Q$. An easy calculation shows that $\mathcal H_1$  has dimension $3.$ 

Let  $\mathcal H_2$ be the closure in $\subset V_{Q,\theta}$ of the locus corresponding to quartics $C$ with a node $p\in Q$. 

Denote $\Pi: {\Ps4}^* \dashrightarrow V$ the morphism associating to a hyperplane section $H$, such that $\ell \not \subset  H$, the ramification curve of the resolution of $\pi_\ell|_{X\cap H}$. We can extend the map $\Pi\Lambda$  to whole $\Ps4$ such that it is the identity map (see Lemma~\ref{qua}).
 
If $(H,D_H)\in  \bar C_{\ell}$ then $\Psi$  is not defined at $(H,D_H)$. Since $\ell\subset H$, the rational map $\pi_{\ell}|_{X_H}$ is not a projection with center in a point, hence we cannot obtain a quartic curve as branch locus.

Let $\tilde{\Pi}: \mathbf{P}' \ra V$ the morphism obtained by blowing up  ${\Ps4}^*$ along the plane parameterizing hyperplanes containing $\ell$.
Denote with ${\tilde C}'=\mathbf{P}' \times_{{\Ps4}^*}\tilde C\stackrel{\nu}{\rightarrow} \tilde C$. Then $\tilde{C}'$ parameterizes triples $(H,D_H,p)$ with $p\in\ell\cap H$.  
The regularized map ${\Psi}':{\tilde C}'\rightarrow V_{Q,\theta}$ maps $(H,D_H,p)$ to the ramification curve of the projection with center $p$ and associates an even theta-characteristic to it. Note that $p$ is on a line of $X\cap H$, hence $C(X \cap H,p)$ is singular.

We denote with $F(X)_{\ell}$ the subvariety of the Fano surface of lines on $X$ given by
\[F(X)_{\ell}:=\{ \tilde{\ell}\in F(X): \ell\cap \tilde{\ell}\ne \emptyset \}.\]
Define
\[\mathcal H'=\{(H,D_H,p)\in {\tilde C}': p\in\ell_H \mbox{ with }\ell_H\in F(X),\ell_H\subset H\},\]  
i.e., the set of hyperplane sections such that $p$ lies on one of the lines of $X$ contained in that hyperplane.
\begin{Prop}\label{H_2} Let $d\in V$ such that $X_d$ is a smooth hyperplane section of $X$, the line $\ell$ is not contained in $X_d$ and the point $p\in \ell \cap X_d$ lies on one of the 27 lines contained in $X_d$. Then the quartic $C_d$ is nodal with at most three nodes. The singularities lie on the common intersection of $C_d$, the quintic $Q$ and the bitangent $b_d$. 

The locus $\mathcal H_2\subset V_{Q,\theta}$  is of dimension $3$.
\end{Prop} 
\begin{proof} Let $\tilde{X}_d$ be the blow-up of $X_d$ in $p.$ Then $\tilde{X}_d$ is a Del Pezzo surface containing at least one and at most three disjoint $-2$-curves, namely the strict transforms of lines through $p$. 

The morphism $\tilde{\pi}_p$ contracts the $-2$-curves to nodes of $C_d$ (see Example~\ref{node}). 
Consider now the projection from $X$ with center $\ell$. Let $\ell'$ be a line intersecting $\ell$. 
 All the points on $\ell'-(\ell'\cap \ell)$ are mapped to one point $q'$ (namely the image of the 2-plane spanned by $\ell$ and $\ell'$). The fiber of $\tilde{\pi}_\ell$ over $q'$ is a reducible conic, hence $q'\in Q$. Since $\ell'$ is mapped to $q'$ we obtain that the node $q$ of $C$ coincides with $q'$. The line $b_d$ parameterizes lines that are tangent at $X_d$ in $p$, so $q'\in b_d$.

The locus $\mathcal H_2$ is birational to a covering of degree $36$ onto the set of hyperplanes containing a line in $F(X)_\ell$. Since $F(X)_\ell$ is a curve in the Fano scheme $F(X)$ (\cite[Lecture 2.1]{Tyu}), it follows that $\mathcal H_2$ is 3-dimensional.
\end{proof}

\begin{Rem}\label{Eck} 
Suppose $\ell_1, \ell_2$ and  $\ell_3$ are coplanar lines having a common intersection point $p$. (Then $p$ is called a {\em star point} or {\em Eckardt point}). If we project from $p$ we obtain a quartic curve $C_d$ with three collinear nodes, hence $C_d$ is reducible and contains the line $b_d$.
Such a $d$ exists, because the locus of cubic surface with a star point in $\{X_d\}_{d \in V}$ has codimension at most 1, so the locus of cubic surfaces $X_d$ with $X_d\cap \ell$ as star point has codimension at most 3.
\end{Rem}

\begin{Lem}\label{noniso} Suppose $H$ is a hyperplane such that $X\cap H$ has non-isolated singularities. Then $X\cap H$ is an irreducible cubic surface containing a double line. Resolving the projection $\tilde{\pi}_p$ from a point $p\in X\cap H$ not on the double line gives rise to a double cover of $\Ps2$ ramified along a reduced conic and a double line, or a triple line union with a line.\end{Lem}

\begin{proof} If $X\cap H$ is reducible then $X$ contains a 2-plane. An easy calculation shows that then $X$ is singular, which contradicts our assumptions on $X$. Hence $X\cap H$ is irreducible. By the classification of cubic surfaces it follows that the singular locus of $X\cap H$ is a double line. 

Since $X_H$ can be obtained in a family of smooth cubic surfaces, the ramification locus of $\tilde{\pi}_p$ is a quartic curve $C$.
Since the singular locus of the double cover is mapped to the singular locus of $C$, we obtain that $C$ is non-reduced. If $C$ were a double conic or a line with multiplicity four, then the double cover would be reducible. So $C$ is either the union of a conic and a double line or the union of a line and a triple line. 
\end{proof}

We describe what type of non-reduced quartics occur in $\{C_d\}_{d\in V}$. 
Let \[N:=\{d\in V: C_d \mbox{ is non-reduced}\}.\]
\begin{Lem}\label{dc}
The family $V$ contains no $d$ such that $C_d$ is a double (possibly singular) conic.
\end{Lem}
\begin{proof} Suppose $d\in |K_Q \otimes \eta|$ is such that $C_d$ is a double conic. This means that $2d$ is cut out by the double conic, hence $d$ is cut out by a conic, so $d \in |\str_Q(2)|=|K_Q|$. This implies that $\eta \sim 0$, hence $\theta \sim \str(1)$. This cannot happen, since $h^0(\theta)=1$.
\end{proof}
Let $T$ be the conic cutting out on $Q$ the effective divisor of $\theta$.  
\begin{Prop}\label{PC} 
The morphism $\Lambda$ maps $N$ isomorphically to the set of hyperplanes in $\Ps4$ containing the line $\ell$.
In fact, $N=\{\theta+t: t\in \mid \str(1)\mid\}$.
The quartic curves $C_d$ with $d\in N$ are the union of the conic $T$ and a double line.
\end{Prop} 
\begin{proof} 
 From Lemma~\ref{dc} it follows that $d\in N$ if and only if the quartic $C_d$ contains at least a double line $r$. This means that $d-(r.Q)\sim \theta$ is effective,  $d=(r.Q)+(T.Q)$ and $C_d=2r+T$. Thus $N=\theta+\mid \str(1)\mid\cong\Ps2$ and the last assertion follows.

Let $\mathcal I$ be the 2-dimensional subvariety in ${\Ps4}^*$ of hyperplanes containing $\ell$. Let $H\in \mathcal I$ and $\Xi:=Q'\cap \ell$ be the set of vertexes on the line $\ell$. We have seen (cf. Lemma~\ref{con}) that the conic $T$ parameterizes plane sections $P=C_p\cup \ell$ such that the conic $C_p$ is tangent to $\ell$. Then the set $\Xi$ parameterizes plane sections $P$ over $T\cap Q$.  Let $d\in V$ be such that $\Lambda(d)=H$. If the support of $d$ contains $T\cap Q$, then $d\in \theta+\mid \str(1)\mid$ and $C_d$ is the union of $T$ and a double line $r$. The pencil of planes in $H$ containing $\ell$ gives a line $h=\pi_{\ell}(H)$ in $\Ps2$. The intersection $Q\cap h$ is the projection of vertexes in $H\backslash\ell$, so $r=h$.  

Conversely, let $d\in N$ such that $2d$ is cut out by the non-reduced quartic $C_d=T+2r$. 
The corresponding hyperplane $H_d$ contains the 5 vertexes in $Q'\cap \ell$, hence $H$ contains $\ell$. In fact, $H_d$ is the unique hyperplane containing $\ell$ and the pencil of planes corresponding to $r$.
\end{proof} 

\begin{Rem}
The conic $T$ can be reducible. In this case there are $d_1,d_2\in N$  such that $C_{d_1}$ and $C_{d_2}$ are the union of a triple line and a line.

The five points in $Q\cap T$ correspond to plane sections of $X$ with a vertex on $\ell$. These are the  Eckardt planes in $X$ containing $\ell$.
\end{Rem}

\begin{Prop}\label{strangecon} There does not exist an $H\in {\Ps4}^*$ such that $X \cap H$ has non-isolated singularities. Moreover, if $d\in \Ima(\Pi)$ then $C_d$ is a quartic curve with at most $ADE$ singularities. \end{Prop}

\begin{proof} If $X \cap H$ had non-isolated singularities, then by Lemma~\ref{noniso}, there would be a $d\in \Ima(\Pi)$ with $C_d$ a conic with a double line or a triple line with a line. However, by Proposition~\ref{PC}, non-reduced quartics in $V$ are not in the image of $\Pi$. 
\end{proof}
Consider the diagram:
\[\begin{array}{ccccc}
 {\tilde C}' & \stackrel{\nu}{\longrightarrow} & \tilde C & \stackrel{\tilde{\mathcal A}}{\longrightarrow} & \tilde{\Theta} \\
\downarrow {\Psi}' &&&& \tilde{\G}\downarrow \\
V_{Q,\theta}&\stackrel{\mathcal F}{\longrightarrow}& V &\stackrel{\Lambda}{\longrightarrow} &{\Ps4}^* . \end{array}  \]
We define
\[ X'_F:=\{H\in {\Ps4}^*: \exists \ell' \in F(X)_\ell, \ell' \subset H \}\]

\begin{Prop} 
The following properties hold:
\begin{enumerate}
\item the branch locus of $\tilde{\pi}=\tilde {\mathcal G}\tilde{\mathcal A } $ is $X'$ and  $\tilde{\pi}^{-1}(X')=\mathcal H;$ 
\item the map ${\Psi}'$ induces an isomorphism  from ${\tilde C}'\backslash (\mathcal H')$ to $V_{Q,\theta}\backslash\mathcal H_2$;
\item the locus  in $V$ of points $p$ such that  $\mathcal F^{-1}(p)$ has less then 36 elements is $\Lambda^{-1}(X'\cup X'_F)$. Moreover,  $(\Lambda\mathcal F)^{-1}(X'\cup X'_F)=\mathcal H_1\cup \mathcal H_2.$
\end{enumerate}
\end{Prop}

\begin{proof} 

\begin{enumerate}
\item The map $\tilde{\pi}$ is regular on the open subset corresponding to nodal hyperplane sections (see~\ref{model}), where it is just the projection $(H,D_H)\ra H$. Then the result follows from Example~\ref{node}.
\item This follows from Proposition \ref{doublesix}.
\item The map $V_{Q,\theta}\ra V$ is ramified over the locus of $d\in V$ such that $C_d$ is singular. These $d$ correspond to hyperplanes $H$ such that either $X_H$ is singular or $H$ contains a line in $F(X)_{\ell}$.
\end{enumerate}
\end{proof}

It is natural to consider the stratification of the new model $V_{Q,\theta}$ corresponding to the number of nodes of the quartics.
So we define $ V_{Q,\theta}(\delta)$ the subset of the model whose general element has a support curve with $\delta$ nodes. In this way we find a stratification induced in the two divisors corresponding to singular quartics. We denote the induced strata with $\mathcal H_i(\delta):=\mathcal H_i \cap V_{Q,\theta}(\delta).$ 
\begin{Prop}
The subset of the quartics in $V$ which are the union of a cubic $K$ and a line $M$ such that $K\cap M \subset Q$, is isomorphic to the blow-up of $Q$ in 5 points. In particular, the corresponding subset of $\mathcal H_2(3)$ is one dimensional.
\end{Prop} 
\begin{proof}
Consider a line $l$ tangent to $Q$ with $l.Q=p+q+r+2s$. The set of cubics $K$ passing transversally through the points $p,q,r$ with $K+l\in V$ is parameterized by linear systems $L$ with:
\[2L\sim \str (3)-p-q-r\sim \str(2)+2s,\]
\[L+\str(1)-s\sim\theta(1).\]
Then $L\sim \theta+s$. Let $Q\cap T=\{p_1,\dots,p_5\}$ then $h^0(\theta+p_i)=2$ with $i=1,\dots,5$.

If $s\in Q-(Q\cap T)$ then, $h^0(\theta-s)=0$, because $\theta$ has a unique section which is non-zero at $s$, from the Riemann-Roch Theorem it follows that then $h^0(\theta+s)=1$. With the same type of argument it follows that $h^0(\theta+p_i)=2$ for $i=1,\dots 5$.

Consider the morphism $H_2(3)\stackrel{\phi}{\rightarrow}Q$ defined as $\phi(l,K)=(l\cap Q)\backslash (K\cap Q).$
The fiber of $\phi$ over $s$ is isomorphic to $\Ps{}(H^0(\theta+s))$. Then $\phi$ is the blow-up of $Q$ in $p_1,\dots,p_5$, in particular the variety $H_2(3)$ is 1-dimensional.
\end{proof}

\begin{Rem}
The special interest for $\mathcal H_2(3)$ lies in the fact that, as we observed in Remark~\ref{Eck}, the existence of quartics of a cubic and a line in $\mathcal H_2$ corresponds to the existence of hyperplanes $H$ such that the section $X\cap H$ contains an Eckardt point $p$ and $\{p\}=\ell\cap H$.\end{Rem}

As $V_{Q,\theta}$ is a fourfold it is natural to expect that the stratum $V_{Q,\theta}(4)$ is a zero dimensional scheme. We count the number of elements in $V_{Q,\theta}(4)$.

\begin{Def}
A triple $\theta_1,\theta_2,\theta_3$ of odd theta-characteristics on a smooth quintic curve $Q$ is called {\em syzygetic} if and only if $\theta_1+\theta_2-\theta_3$ is an odd theta characteristic.\end{Def}
Note that three odd theta characteristics are syzygetic if and only if their points of contact with $Q$ lie on a quartic curve.
Each theta characteristic $\theta$ on $Q$ defines a quadratic form $q_\theta$ on $\Jac(Q)[2]$:
\[ q_{\theta}(\eta)=h^0(\theta+\eta)+h^0(\theta)\bmod  2. \]
This quadratic form it closely related to the Weil-pairing $\langle \cdot, \cdot \rangle$   on $\Jac(Q)[2]$ by the Riemann-Mumford relation:
\[ q_\theta(\alpha+\beta)+q_\theta(\alpha)+q_\theta(\beta)=\langle \alpha,\beta \rangle, \; \alpha,\beta \in \Jac(Q)[2],\]
see for example \cite[page 290]{ACGH}.
\begin{Lem}\label{syzlem}
A triple of odd theta characteristics $\theta_1,\theta_2,\theta_3$ is syzygetic if and only if:
\[q_{\theta_1}(\theta_2-\theta_3)=0.\]
\end{Lem}
\begin{proof} We have that $q_{\theta_1}(\theta_2-\theta_3)=0$ is equivalent to 
\[ h^0(\theta_1+\theta_2-\theta_3)\equiv h^0(\theta_1) \bmod 2. \]
Since $h^0(\theta_1)\equiv 1 \bmod 2$, the above is equivalent to $h^0(\theta_1+\theta_2-\theta_3)\equiv 1 \bmod 2$.
\end{proof}
\begin{Prop}
There are 495 curves appearing as supports of spin curves in $\mathcal H_1(4)$ corresponding to unordered pairs of conics $\{C_1,C_2\}$ such that $(C_1.Q)_p+(C_2.Q)_p$ is even for all $p\in Q\cap (C_1\cup C_2).$   
\end{Prop}
\begin{proof} Suppose $D\in 2 \Div(Q)$ then we denote $\red(D):=\frac{1}{2}D$. 
 If $C$ is a conic everywhere tangent to $Q,$ then $\red(C\cdot Q)$ is an odd theta characteristic on $Q$. 
We are looking for all unordered  couples $(C_3,C_4)$ of distinct conics such that
$\red(C_3\cdot Q+ C_4 \cdot Q)\in \theta(1)$. Let $\theta_3$ and $\theta_4$ be the associated odd theta-characteristics. Then 
$\theta_3 + \theta_4\sim \theta(1)$.

Let $\theta_1=\theta$ and $\theta_2=\mathcal O(1)$. We are looking for all $\theta_3\neq \theta,\str(1)$ such that $\theta_1,\theta_2,\theta_3$ are syzygetic. Given these three, $\theta_4\sim \theta_1+\theta_2-\theta_3$.

Define $\eta=\theta_1-\theta_2$. Consider  $Z=(q_{\theta_1}^{-1}(0)\cap q_{\theta_2}^{-1}(0))\backslash \{\eta,0\}$. From  Lemma~\ref{syzlem} it follows that the elements in $\alpha \in Z$ give rise to syzygetic triple $(\theta_1,\theta_2,\theta_2+\alpha)$, and every syzygetic triple $(\theta_1,\theta_2,\theta_3)$ gives rise to an element $\theta_3-\theta_2 \in Z$.

We have 
\[ q_{\theta_2}(\alpha)=q_{\theta_1}(\alpha)+\langle \alpha, \eta \rangle\]
for every $\alpha \in \Jac(Q)[2]$ (cf. \cite[Corollary I.3.21]{Ver}).

Let $\alpha\in Z$, then $\langle \alpha, \eta \rangle=0$. Thus $\alpha$ is an element of $V=(\F_{2 }\eta)^{\bot}/\F_{2} \eta\cong {\F_{2 }}^{10}$. We have also $q_{\theta_1}(\eta)=q_{\theta_2}(\eta)=0$. It can be shown that the induced quadratic forms on $V$ coincide and give an odd quadratic form $q$.
The fiber $q^{-1}(0)$ contains $2^{4}(2^5-1)=496$ elements (see for example \cite[Section 1]{GrHa}). We get $2\cdot 495=990$ elements in $Z$. To any (unordered) pair of conics correspond two elements in $Z$, we get 495 pairs of conics in total.
\end{proof}

\section{New model and stable reduction of algebraic curves}\label{stablered}
In this section we explain one consequence of the existence of our model for the theta divisor for the stable reduction of non-reduced quartic curves.

First, we recall the stable reduction theorem, a proof can be found in \cite[proof of Theorem 1.1]{Ba}. In this section we denote with $\Delta$ the disc $\Delta:=\{t\in\C:|t|\le 1\}.$
\begin{Thm}
Let $f:\mathcal C\rightarrow \Delta$ be a  flat family with $C_{t}=f^{-1}(t)$ a Deligne-Mumford stable curve of genus $g\ge 2$ for $0\ne t \in \Delta.$
Then there exists a commutative diagram:
\[\begin{array}{ccccc}
 \mathcal{Y} & \stackrel{\phi}{\longrightarrow} & \mathcal C\times_\Delta \Delta & \rightarrow & \mathcal C \\
\downarrow\tilde{f} &&\downarrow f\times_\Delta p&& \downarrow f \\
\Delta &=&\Delta&\stackrel{p}{\longrightarrow} & \Delta \end{array}  \]
such that
\begin{enumerate}
\item  the morphism $p: \Delta \ra \Delta$ is given by $z \mapsto z^k$ for some integer $k>0$.
\item  the extension $\tilde f$ is a family of Deligne-Mumford stable curves.
\item  $\phi$ is an isomorphism away from the fiber over $0\in \Delta.$
\end{enumerate}
Moreover, any two extensions $(\tilde{f},p)$ and $(\tilde{g},q)$ satisfying the above three conditions have isomorphic fibers over $0$.
\end{Thm} 

\begin{Def} The {\em stable reduction of $C_0$ with respect to $f$} is the curve $\tilde{f}^{-1}(0)$ and denoted by $\mathcal R_{f}(C_0)$.\end{Def} 

Note that the stable reduction depends on the family chosen. If the family $f$ is not a general smoothing of the curve $C_{0},$ it is usually difficult to calculate the stable reduction.

We consider now the case of the union $C_0$ of a conic $T$ and a double line $L.$ 

\begin{Prop}There exist infinitely many flat families $f:\mathcal C\rightarrow \Delta$ of plane quartics whose central fiber is $f^{-1}(0)=C_0$ and $f^{-1}(t)$ is a Deligne-Mumford stable quartic, such that $\mathcal{R}_{f}(C_0)$ is a quartic with one node and everywhere tangent to a suitable quintic $Q$ with the node on $Q.$ 
\end{Prop}
\begin{proof}
Let $Q$ be a smooth  quintic curve $Q$  everywhere tangent to $T$ and satisfying the conditions mentioned at the beginning of Section~\ref{newmodel}.  The divisor $\frac{1}{2}(T\cdot Q)$ gives an odd theta characteristic $\theta$ on $Q.$ Denote with $X$ the corresponding cubic threefold  and $\ell$ the line on it.  
Let $V=\mid \theta(1)\mid, \Pi, \Lambda$ as in Section~\ref{newmodel}. 
Note that $C_0$ belongs to $V$. Let $\mathcal H_{2}$ be the divisor of quartics in $V$ with a singular point on $Q.$

From Proposition~\ref{PC} follows that $\Lambda(C_t)$ is a hyperplane not containing $\ell$ for $t\neq 0$ and containing $\ell$ for $t=0$. The family $\Lambda(C_t)$, $t\in \Delta$ gives a curve $\Sigma$ in ${\Ps4}^*$. Let $\tilde{\Sigma}$ be the strict transform of $\Sigma$ in $\mathbf{P}'$. Let $\Sigma_0$ be the intersection of $\tilde{\Sigma}$ with the exceptional divisor. Then by Proposition~\ref{H_2}, we obtain $\tilde{\Pi}(\Sigma_0)\in \mathcal{H}_2$. Since $\tilde{\Pi}\Lambda(C_t)=C_t$, for $t\neq 0$, we obtain that if $\tilde{\Pi}(\Lambda_0)$ is a Deligne-Mumford stable curve, then $\mathcal R_{f}(C_0)=\tilde{\Pi}(\Sigma_0)$.

If we vary $Q$ and $f$ enough then we can find families such that $\tilde{\Pi}(\Lambda_0)$ is a Deligne-Mumford stable curve. If we vary even more, we obtain that $\tilde{\Pi}(\Lambda_0)$ is a generic element of the associated $\mathcal{H}_2$, hence is a quartic curve with a node on $Q$.
\end{proof}

\end{document}